\documentclass[a4paper,12pt]{article}
\usepackage{t1enc}
\usepackage[latin1]{inputenc}
\usepackage[french]{babel}
\usepackage{amssymb,latexsym,amsfonts}
\usepackage{amsmath,amssymb}
\usepackage{url}

\newtheorem{thm}{Théorème}
\newtheorem{pro}{Proposition}
\newtheorem{cor}{Corollaire}
\newtheorem{lem}{Lemme}

\def \proj{\mathbb{P}}
\def \nat{\mathbb{N}}

\def \ree{\mathbb{R}}
\def \com{\mathbb{C}}

\def \polyr{\ree[x_{1},\ldots,x_{n}]}

\def \eps{\varepsilon}
\def \cinf{\mathcal{C}^{\infty}}

\def \ppoe{\leqslant}
\def \pgoe{\geqslant}

\newcommand{\uplet}[2]{\ensuremath{(#1_{1}, \ldots ,#1_{#2})}}
\newcommand{\iuplet}[3]{\ensuremath{(#1_{#2}, \ldots ,#1_{#3})}}
\newcommand{\polyn}[3]{\ensuremath{#1 \lbrack #2_{1}, \ldots , #2_{#3}\rbrack }}
\newcommand{\norm}[1]{\ensuremath{ \left \Vert #1\right \Vert}}
\newcommand{\modul}[1]{\ensuremath{ \left \vert #1\right \vert}}

\newcommand{\prim}[1]{\ensuremath{#1^{\prime}}}
\newcommand{\prims}[1]{\ensuremath{#1^{\prime \prime}}}
\newcommand{\derive}[2]{\ensuremath{\frac{\partial}{\partial #2}\! #1}}

\newcommand{\braces}[1]{\ensuremath{\left \lbrace #1 \right \rbrace}}

\newcommand{\parentheses}[1]{\ensuremath{\left ( #1 \right )}}

\newcommand{\restric}[2]{\ensuremath{#1\vert_{#2}}}
\newcommand{\bracks}[1]{\ensuremath{\left \lbrack #1 \right \rbrack}}
\newcommand{\obracks}[1]{\ensuremath{\left \rbrack #1 \right \lbrack}}
\newcommand{\cobracks}[1]{\ensuremath{\left \lbrack #1 \right \lbrack}}

\begin{document}
\title{\Large Inégalités de Markov tangentielles locales sur les courbes algébriques
singulières de $\ree^{n}$\\   \normalsize \textit{version
révisée}\footnote{Prépublication n°109 (1997) Laboratoire E.PICARD C.N.R.S. U.M.R. 5580
}\\ \Large par}
\author{\textsc{Laurent GENDRE}}
\date{ }
\maketitle

\textbf{Abstract-} We prove that all real singular algebraic curves admits Markov 's
local tangential inequalities. We give a geometric significance of Markov's exponent.
\section{Introduction}
  Ce sont dans les articles de L.BOS, P.MILMAN, N.LEVENBERG et B.A. TAYLOR \cite{bos2} et
 W. PLESNIAK et M. BARAN \cite{baran1} qu'il est montré que l'exposant de Markov est au moins
  égal à 1 sur les surfaces algébriques de $\ree^{n}$, lisses $\mathcal{C^{\infty}}$, sans bord et compactes.
  Pour y parvenir, les auteurs utilisent essentiellement des outils de la théorie du pluripotentiel
complexe. En revanche dans l'article R. NARASIMHAN et C. FEFFERMAN \cite{fefferman}, il
est  démontré qu'il existe de telles inégalités avec un exposant de Markov valant 2 dans
des compacts, ne contenant pas de points singuliers, inclus dans des sous-ensembles
algébriques de $\com^{n}$, et ceci par des techniques non triviales de géométrie
analytique. Bien sûr, dans cette introduction trop courte, on ne pourra citer toutes les
parutions et  les auteurs ayant contribué à l'essor de ce sujet vu la densité de
résultats. Ceci dit, donnons une liste non exhaustive des résultats intéressants et
fondamentaux de ce domaine qui permettra de cerner l'ensemble :\cite{pawlucki1},
\cite{pawlucki2}, \cite{baran2}, \cite{baran3}, \cite{bos2}, \cite{bos3}, \cite{brudnyi}
et \cite{zeriahi}.

Dans cet article, nous démontrons l'existence d'inégalités de Markov locales
tangentielles sur les courbes algébriques singulières de $\ree^{n}$. Nous répondons sur
le comportement de l'exposant de Markov, en montrant qu'il est minoré par un invariant
géométrique. Pour cela, nous utilisons des outils de géométrie analytique et des
résultats de la théorie du pluripotentiel complexe,  afin d'affaiblir la notion d'$HCP$
pour la fonction de Green avec pôle à l'infini.

\section{Énoncé des Théorèmes}
 Nous adopterons l'identification suivante
 $\com^{n}\equiv \ree^{n}\oplus i \ree^{n}$; $\ree^{n}$ est donc un sous-espace vectoriel de
 $\com^{n }$. Par conséquent, nous considérons, canoniquement, $\ree[x_{1},\ldots,x_{n}]$
 comme sous-espace de $\com[z_{1},\ldots,z_{n}]$.

  Si $E$ est un sous-ensemble de $\com^{n}$ et $a$ un point de
$\overline{E}$, on dira que $v$ est un vecteur tangent en $a$ de $E$, s'il existe une
suite $(a_{i})_{j\in\nat}$ de point de $E$ et une suite de réels strictement positifs
$(\lambda_{j})_{j\in\nat}$ telles que $v=\lim_{j\rightarrow+\infty}\lambda_{j}(a-a_{j})$.
Nous appellerons l'ensemble des vecteurs tangents de $E$ en $a$ le cône tangent en de $E$
en $a$ et le noterons $C(E,a)$. Notons que pour tout $v\in C(E,a)$, $\lambda\cdot v\in
C(E,a),\quad\forall\lambda\in\ree_*^+$;  les vecteurs tangents \textit{unitaires} de
$C(E,a)$ ont un sens.

Pour tout sous-ensemble $K\subset\subset\com^{n}$, nous noterons la norme uniforme sur
$K$ ainsi: $\norm{.}_{K}$.

Nous appellerons \textit{morceau de courbe algébrique}, tout sous-ensemble analytique
connexe et irréductible d'une courbe algébrique. Nous conviendrons de la notation
suivante : $\com_{i}=\com^{i}\times \lbrace 0 \rbrace$ sous-espace de $\com^{i} \times
\com^{n-i}$. Enfin pour tout sous-ensemble algébrique $A$ de $\ree^{n}$, nous noterons
respectivement l'ensemble des points singuliers et réguliers par $A_{sing}$ et $A_{reg}$.
\begin{thm}
Soit $A$ un morceau de courbe algébrique. Pour tout $x_{0}$ dans $A$, il existe
$C_{1},C_{2}$ et $\eps_{0}$ des constantes strictement positives, dépendant de $x_{0}$ et
localement majorées telles que $\forall \eps \in \lbrack 0,\eps_{0}\lbrack,\;\forall p
\in\polyn{\ree}{x}{n}$:
\begin{enumerate}
\item $\;si\;x_{0}\in A_{sing}$,

 $$
     \vert D_{v}p(x_{0}) \vert  \leqslant C_{1} \left ( \frac{C_{2} (\deg(p))^{2}}{\eps} \right
     )^{k} \norm{p}_{A \cap B(x_{0},\eps^{k})},
 $$
où $k$ est la multiplicité complexe du point singulier $x_{0}$ dans $\tilde{A}$ et v un
vecteur unitaire dans $C(A,x_{0})$.

\item $\;si\;x_{0}\in A_{reg}\setminus \partial A,$

 $$
      \vert D_{v}p(x_{0}) \vert  \leqslant C_{1}
      \left ( \frac{C_{2} \deg(p)}{\eps} \right )^{}
      \norm{p}_{A \cap B(x_{0},\eps^{})},
 $$

\item $\;si\;x_{0}\in \partial A\setminus A_{sing},$

 $$
       \vert D_{v}p(x_{0}) \vert  \leqslant C_{1}
      \left ( \frac{C_{2} \deg(p)}{\eps} \right )^{2}
       \norm{p}_{A \cap B(x_{0},\eps^{})},
 $$
 où $v$ est un vecteur unitaire de l'espace tangent  $T_{x_{0}}A_{reg}$.
 \end{enumerate}
\end{thm}
\begin{cor}
Soit $A$ une courbe algébrique, compact, localement irréductible et sans singularité au
bord. Alors A admet des inégalités de Markov tangentielles locales.  Il existe donc
$C_{1},C_{2}$ et $\eps_{0}$ des constantes absolues strictement positives telles que :
 $\forall \eps \in \lbrack 0,\eps_{0}\lbrack ,\; \forall x_{0} \in A, \; \forall p \in \polyr,$
$$
\vert D_{v}p(x_{0}) \vert  \leqslant C_{1} \left ( \frac{C_{2} (\deg(p))^{2}}{\eps}
\right )^{ k}  \Vert p \Vert_{A \cap B(x_{0},\eps^{k})},
$$
 où $k$ est la multiplicité
complexe du point $x_{0}$ et $v$ un vecteur unitaire dans $C(A,x_{0})$.
\end{cor}

Au 2 du théorème 1, nous retrouvons le résultat  énoncé par Fefferman et Narasimhan dans
\cite{fefferman}, pour le cas d'une courbe algébrique.
 Les hypothèses du théorème 1 englobent, bien sûr, les courbes algébriques lisses sans bord. Nous en déduisons
 le corollaire suivant qui est une conséquence du théorème de Bos, Milman, Levenberg et Taylor \cite{bos2}.
\begin{cor}
Soit $A$ une courbe algébrique $\cinf$, compacte et localement irréductible de
$\ree^{n}$. $A$ admet des inégalités de Markov tangentielles locales d'exposant 1.
Autrement dit, il existe $C_{1},C_{2}$ et $ \eps_{0}$ des constantes positives absolues
telles que : $\forall \eps \in \lbrack 0,\eps_{0}\lbrack,\; \forall x_{0} \in A,\;
\forall p \in \polyn{\ree}{x}{n}$,
     $$
        \vert D_{v}p(x_{0}) \vert \leqslant C_{1} \left ( \frac{C_{2} \deg(p)}{\eps}
        \right ) \Vert p \Vert_{A \cap B(x_{0},\eps)},
     $$
 où $v$ est un vecteur unitaire de l'espace tangent $T_{x_{0}}A$.
\end{cor}
 \section{Préliminaires}

 \subsection{Sous-ensembles algébriques de $\ree^{n}$}

 Si $S$ est une partie de $K^{n}(n\pgoe 2)$ ($K$ est le corps $\ree$ ou $\com$), on notera que
 $$
 I^K(S):=\braces{p\in\polyn{K}{x}{n}:\restric{p}{S}=0},
 $$
 est un idéal de $K^{n}$ ayant un nombre fini de générateur ($\polyn{K}{x}{n}$ est un anneau
 Noetherien). Si $\mathcal{P}$ est une partie non vide $\polyn{K}{x}{n}$, on écrira
 $$
    loc\;\mathcal{P}=\braces{x\in K^{n}:\forall p\in\mathcal{P},\;p(x)=0}
 $$
le locus de $\mathcal{P}$.

 Les sous-ensembles algébriques de $K^{n}$ sont les parties $A$
de $K^{n}$ telles que:
$$
       loc\;I^K(A)=A.
$$

  En gardant l'identification  $\com^{n}\equiv\ree^{n}\oplus i\ree^{n}$, nous
  considérerons, pour tout sous-ensemble algébrique $A$ de $\ree^{n}$ le complexifié de $A$ comme étant le
  plus petit sous-ensemble algébrique complexe de $\com^{n}$ contenant $A$. On le notera
  $\tilde{A}$, ainsi que, $\tilde{A}_{reg}$ et $\tilde{A}_{sing}$ les
  sous-ensembles, respectifs, singuliers et réguliers de $\tilde{A}$. On remarquera que $I^{\com}(\tilde{A})=I^{\ree}(A)\otimes_{\ree}\com$ et
  $A=\tilde{A}\cap\ree^{n}$. Pour plus
  de précisions, il faut se référer au livre de Narasimhan(\cite{narasimhan},
  page 91).

  Un sous-ensemble $A$ de $\ree^{n}$ est une courbe algébrique, s'il existe des polynômes
  $p_{1},\ldots,p_{s}$ dans $\polyn{\ree}{x}{n}$, ($s\in \nat^{\star}$), tels que :
  $A= \lbrace x \in \ree^{n} :p_{1}(x)=\cdots=p_{s}(x)=0 \rbrace$ et que $\dim_{\ree} A=1$.
 \subsection{Critère d'algébricité}

Le critère d'algébricité de Rudin sera utile dans  la démonstration du théorème 1. Une
seule implication nous servira. On trouvera sa démonstration dans(\cite{chirka}, Théorème
3, page 78), voici son énoncé :
\begin{thm}
Soit $X$ un sous-ensemble analytique de $\com^{n}$($n \in \nat $) de dimension pure $p$.
Alors $X$ est algébrique si et seulement s'il existe $C$ et $s$ deux constantes
strictement positives et un changement de coordonnées unitaire tels que : $$ \forall
z=(z^{\prime},z^{\prime \prime}) \quad \vert z^{\prime \prime} \vert  \leqslant C(1+\vert
z^{\prime} \vert^{s}) $$ où $$ z^{\prime}=(z_{1},\cdots ,z_{p}) \; et \; z^{\prime
\prime}=(z_{p+1},\cdots ,z_{n}) $$
\end{thm}

Nous verrons plus tard que nous choisirons un changement de coordonnées unitaire bien
spécifique et que ce choix est crucial dans l'obtention de nos estimations. En effet,
nous aurons besoin que le sous-espace $\ree^{n}$ de $\com^{n}$ soit stable par ce
changement de coordonnées. En fait, choisir un cône ,comme au théorème ci-dessus, c 'est
trouver une "bonne" projection. Ceci peut être énoncé par la proposition suivante.
\begin{pro}
Soit $\tilde{A}$ un sous-ensemble algébrique de $\com^{n}$ de dimension pure p. Si la
projection $\pi \, : \,\tilde{A}\longrightarrow \com_{p}$ est propre on a : $$ \tilde{A}
\subset \lbrace (\xi^{\prime},\xi^{\prime \prime}) \in \com^{n} : \vert \xi^{\prime
\prime} \vert  \leqslant C(1+\vert \xi^{\prime} \vert^{s})\rbrace. $$ Avec $$
\xi^{\prime}=\uplet \xi p \; et \; \xi^{\prime \prime} = \iuplet \xi {p+1} n $$ où $C$ et
$s$ sont des constantes strictement positives dépendant seulement de $\tilde{A}$.
\end{pro}
Nous nous bornerons juste à énoncer ce théorème dont la démonstration figure dans
(\cite{loja}, Proposition 1, page 389).

\subsection{Multiplicité d'un point dans les sous-ensembles analytiques}

La multiplicité dun point dans un sous-ensemble algébrique complexe de dimension pure $p$
est donné par l'expression suivante :
    $$
           \forall a \in \tilde{A}, \quad \mu_{a}(\tilde{A}):=\min \braces{
           \mu_{a} ( \pi_{L}\vert_{ \tilde{A}}):L \in Gr_{\com}(n, n-p) }.
    $$
où $Gr_{\com}(n,n-p)$ est la Grassmanienne dans $\com^{n}$ et $\pi_{L}$ est la projection
orthogonale définie par $\pi_{L} : \com^{n} \simeq L \oplus L^{\bot} \longrightarrow
L^{\bot}$. On a toujours $\mu_{a}(\tilde{A})\pgoe 1$, pour tout $a\in\tilde{A}$. La
multiplicité complexe $\mu_{a}(\tilde{A})$ n'est autre que le nombre de Lelong du courent
$[\tilde{A}]$ en $a$; c'est donc un invariant par biholomorphisme, cf (\cite{chirka},
 Proposition 2, page 120 et page 190).

\subsection{Paramétrisation de Puiseux pour une courbe analytique complexe}

 Construisons la paramétrisation de Puiseux pour les courbes analytiques de $\com^{n}$.
Nous adapterons la démonstration de (\cite{chirka}, page 67).

Dans cette partie $A$ est un sous-ensemble analytique de $\ree^{n}$ de dimension pure 1,
localement irréductible et $\tilde{A}$ le sous-ensemble algébrique complexifié de $A$
dans $\com^{n}$, tel que $0 \in A_{sing}$ (donc $0\in\tilde{A}_{sing}$); la projection
$\pi:\tilde{A}\cap U\longrightarrow \prim{U}\subset\com_{1}$
($U=\prim{U}\times\prims{U}\subset\com\times\com^{n-1}$ un polydisque) est propre.

Pour toute fonction holomorphe $f$ au voisinage d'un point $a$ de $\com$ nous écrirons
l'ordre de $f$ en $a$: $ord_a\;f$.

Nous noterons l'ensemble $\mathcal{R}^{<r>}$, pour tout suite d'entiers $0\ppoe
l_1<\cdots< l_{r}\ppoe k-1$
 $k\in\nat^{\star}$ et $r\in\nat^{\star}$ tel que  $1\ppoe r \ppoe k$, ainsi :
 $$
 \mathcal{R}^{<r>}=\bigcup_{1\ppoe j\ppoe r}\bracks{0,e^{i\frac{2\pi l_j}{k}}}
 $$
\begin{pro}
Sous les hypothèses et les notations décrites ci-dessus, il existe
$\varphi:\overline{D(0,1)}\longrightarrow \tilde{A}\cap U$ une paramétrisation de Puiseux
définie dans un voisinage ouvert du disque unité fermé $\overline{D(0,1)}$ de $\com_1$
telle que
 $\varphi(z)=(cz^{k},\psi_1(z),\ldots,\psi_n(z)),
\quad\forall z\in \overline{D(0,1)}$, où $k=\mu_{0}(\tilde{A})$, les $\psi_j$ sont
holomorphes avec $ord_{0}\;\psi_j>k$ et $c\in\com$ constante dépendant de $\varphi$; il
existe des entiers $0\ppoe l_1<\cdots< l_{r_A^+}\ppoe k-1$ (resp. $0\ppoe l'_1<\cdots<
l'_{r_A^-}\ppoe k-1$), où $r_A^+$ (resp. $r_A^-$ ) est le nombre de branche de la courbe
réelle $A$ dans $U$, ayant $0$ comme extrémité, au dessus de
$\bracks{0,1}\subset\overline{D(0,1)}$ (resp. au dessus de
$\bracks{-1,0}\subset\overline{D(0,1)}$), de sorte que
$\restric{\varphi}{\mathcal{R}^{<r_{A}^+>}}$ (resp.
$\restric{\varphi}{\mathcal{R}^{<r_{A}^->}}$ ) paramétrise $A\cap
U\cap\braces{\uplet{z}{n}\in\com^{n}:\Re e(z_1)\in\bracks{0,1}}$ (resp. $A\cap
U\cap\\\braces{\uplet{z}{n}\in\com^{n}:\Re e(z_1)\in\bracks{-1,0}})$

\end{pro}

\textbf{Preuve:} $\tilde{A}$  est de dimension complexe 1,
 $\tilde{A}_{sing}$ est un ensemble discret. Donc il existe un polydisque $U=U^{\prime} \times U^{\prime \prime} \subset
 \com \times \com^{n-1}$
tel que $\pi^{<-1>}(0) \cap \tilde{A} \cap U=\lbrace 0 \rbrace$, où $\pi$ est la
projection Nous avons supposé  que le sous-ensemble algébrique $A$ est localement
irréductible; il s'en suit donc que le complexifié $\tilde{A}$ est localement
irréductible, ceci d'après (\cite{narasimhan}, Proposition 2, page 92). La projection
$\pi\vert_{U}$ est propre, d'après le théorème de structure des sous-ensembles
analytiques complexes $(\tilde{A} \cap U,\pi,U^{\prime})$ est un $l$-revêtement
holomorphe ramifié.
 Fixons un point $a \in A \cap U$ tel que $a_{1} \in \ree$ où $a_{1}=\pi(a)$. Soit $\gamma  :\lbrack0,r_{1} \lbrack
\rightarrow \tilde{A}$ le relèvement du segment $\lbrack 0,r_{1} \lbrack \subset
\com_{1}$ passant par $a$ dans $A$ (\textit{i.e} $\exists t_{0} \in \lbrack 0, r_{1}
\lbrack , \; \gamma (t_{0})=a, \;
 \forall t  \in \lbrack 0, r_{1} \lbrack, \; \pi \circ \gamma (t)=t $). Ce relèvement existe car  $(\tilde{A} \cap U,\pi,U^{\prime})$
 est un $l$-revêtement holomorphe ramifié. Comme $\tilde{A}\cap U$ est irréductible $\tilde{A}\cap U \setminus \lbrace 0 \rbrace$
est connexe. Ceci nous permet de dire que $\Gamma_{r}=\pi^{<-1>}(\lbrace z \in \com
:\vert z_{1} \vert =r \rbrace ) \cap \tilde{A}, \; \forall r \in \lbrack 0,r_{1}
\lbrack,$ est une courbe de Jordan fermée. Donc le triplet $(\Gamma_{r},\pi,\vert z_{1}
\vert=r)$ est un $l$-revêtement. On peut ainsi construire une unique paramétrisation
$\gamma_{r}: \lbrack 0, 2 \! \pi \lbrack \rightarrow \Gamma_{r}$ telle que $\pi \circ
\gamma_{r}(t)=re^{i k t}$, avec $\gamma_{r}(0)=\gamma(r).$ Définissons donc:
 $$
           \zeta(z):=\left ( \frac{\pi(z)}{r_{1}} \right )^ {\frac{1}{k}}, \quad
           \zeta(\gamma_{r}(t))= \left ( \frac{r}{r_{1}} \right )^{\frac{1}{k}} e^{it}.
 $$
 La fonction $\zeta$
est holomorphe et injective sur $(\tilde{A} \cap U) \setminus \lbrace 0 \rbrace$. Donc
d'après le théorème de prolongement de Riemann des fonctions holomorphes sur les
singularités, la fonction $z : D(0,1) \rightarrow \tilde{A} \cap U$, telle que $\zeta
\mapsto z(\zeta)$ est une injection holomorphe. Ainsi nous obtenons une paramétrisation
de Puiseux pour $0 \in \tilde{A}_{sing}$.
 Comme le point $a$ n'est pas dans $\pi(A_{sing})$, on a donc
 $\braces{a^{0},\ldots,a^{k-1}}=\pi^{<-1>}(\braces{a_{1}})\cap \tilde{A}\cap U$, en prenant par exemple
  $a=a^{0}$. D'après la construction de la fonction $z(\zeta)$, nous pouvons ordonner les points
$(a^{l})_{l\in\braces{0,\ldots,k-1}}$, de sorte que
 $\zeta(a^{l})\in \lbrack 0,e^{i\frac{2i\pi l}{k}}\rbrack,\;\forall
 l\in\braces{0,\ldots,k-1}$. Définissons les ensembles ci-dessous:
    $$
    \mathcal{E}^{<k>}:=\bigcup_{0\ppoe l \ppoe k-1}\left\lbrack 0,e^{i\frac{2\pi
    l}{k}}\right\rbrack.
    $$
    et $\mathcal{R}^{<k>}$ la partie de $\mathcal{E}^{<k>}$ telle que
    \begin{equation}
       \mathcal{R}^{<r_A^+>}:=\bigcup_{1\ppoe j \ppoe r_{A}^+}\left\lbrack 0,e^{i\frac{2\pi
    l_{j}}{k}}\right\rbrack,
    \end{equation}
où $r_{A}^+$ est le nombre de branches, ayant pour extrémité $0$, dans  $A\cap U$ et
$[0,e^{i\frac{2\pi l_{j}}{k}}]$ est un segment paramétrisant une branche de $A\cap U$
au-dessus du segment $[0,1]$. Bien sûr, l'entier naturel $r_{A}^+$ ne peut dépasser $k$
et est toujours plus grand que $1$. Ainsi $\varphi\vert_{\mathcal{R}^{<r_A^+>}}$
paramétrise analytiquement la partie de la courbe réelle $A$ au-dessus du segment
$[0,1]$. Pour les branches de $A\cap U$ au-dessus de  $[-1,0]$, on se
 ramène au cas des branches au dessus de  $[0,1]$ par une rotation d'angle $\pi$. $\blacksquare$

  Si $v$ est un vecteur
 tangent dans $C(A,0)$. Toutes les branches d'extrémité $0$ de la courbe $A\cap U$ sont paramétrisées localement en $0$ par
 $\varphi\vert_{\mathcal{R}^{<r_A^+>}}$ ou $\varphi\vert_{\mathcal{R}^{<r_A^->}}$.
 Donc pour le vecteur tangent $v$, il correspondra une branche de $A\cap U$,
 à  laquelle il est tangent, et donc un segment de $\mathcal{R}^{<r_A^+>}$ et
 $\mathcal{R}^{<r_A^->}$, du type $[0,e^{i\frac{2\pi l}{k}}]$
  paramétrisant la dite branche par $\restric{\varphi}{\bracks{0,e^{i\frac{2\pi l}{k}}}}$. Bien sûr,
  la géométrie de $\mathcal{R}^{<r_A^+>}$ et de $\mathcal{R}^{<r_A^->}$ dépendent de la
 nature de la singularité de la courbe algébrique réelle $A$.

\begin{lem}
Soit $\tilde{A}$ une courbe  algébrique de $\com^{n}$ irréductible et $a \in
\tilde{A}_{sing}$. Si $k$ est l'exposant de la paramétrisation de Puiseux, alors $k$ est
aussi la multiplicité de la singularité complexe de $\tilde{A}$ en $a$.
\end{lem}

\textbf{Preuve:}Définissons le produit scalaire usuel dans $\com^{n} : \forall (u,v) \in
(\com^{n})^{2}$ avec $u=\uplet{u}{n}$ et $v=\uplet{v}{n}$
    $$
            (u \vert v):= \sum_{j=1}^{j=n} u_{j}
            \overline{v}_{j}.
    $$
 Soient $L \in Gr_{\com}(n-1,n), \; u \in L^{\bot}$ et $v^{2}, \cdots , v^{n}
\in L$ tel que $(u,v^{2}, \ldots , v^{n})$ soit une base $\com$-orthonormée. On a bien:
    $$
            \com u\overset{\bot}{\oplus} \left ( \bigoplus_{2\leqslant
              j \leqslant n}^{\bot} \com v^{j} \right ) =\com^{n},
    $$
on en déduit donc que:
    $$
           \forall z \in \com^{n}, \quad \pi_{L}(z)=(z \vert u) u, \qquad
           \pi_{L}:L \oplus L^{\bot} \longrightarrow L^{\bot}.
    $$
 On peut supposer que  $a=0$, sans perdre de
généralité. Soit la paramétrisation de Puiseux de $\tilde{A}$ en $0$, $z(\xi)=(c
\xi^{k},\psi _{2}(\xi), \ldots , \psi_{n}(\xi))$,
 où les $\psi _{j}$ sont holomorphes dans $D(0,1)$ telles que $ord_{0}(\psi_{j})> k, \; \;
  \forall j \in \braces{2,\ldots, n}$.
Dans (\cite{chirka}, Lemme 1, page 107), nous avons l'égalité suivante:
    $$
          ord_{0}(\pi_{L} \circ
          z(\xi))=\mu_{0}(\pi_{L}\vert _{\tilde{A}}). $$ En identifiant $\com u$ à $\com$, on a $$ \pi \circ
          z(\xi)= \left ( c \! \xi^{k}u_{1} + \sum_{j=2}^{j=n}\psi_{j}(\xi)\overline{u}_{j} \right ) u, \quad
          \forall L \in Gr_{\com}(n-1,n).
    $$
On déduit trivialement que $\mu_{0}(\tilde{A})=k.\;\blacksquare$
\subsection{Fonction de Green avec pôle à l'infini, inégalité de
Bernstein-Walsh et compact HCP de $\com^{n}$}

Soit $E$ un sous-ensemble de $\com^{n}$. On notera  la classe de Lelong,
$L_{E}(\com^{n})$, que l'on définira comme ci-dessous:
      $$
           L_{E}(\com^{n}):=
           \braces{ u \in PSH(\com^{n}): u\vert_{E}\leqslant 0,\; \exists c_{u} \in \ree , \quad u \leqslant
           c_{u} + \log ( 1 + \vert z \vert) },
      $$
     où $ \modul{z}: = \max_{1\leqslant  j\leqslant n} \modul {z_{j}}$, avec $z=\uplet{z}{n}$.

      La fonction de Green avec pôle la l'infini associée  au compact $E$ est  définie ainsi:
      $$
               V_{E}(z):= \sup \braces{ u(z): u\in L_{E}(\com^{n})}, \quad \forall z \in \com^{n}.
      $$
      Il est possible de réduire l'ensemble $L_{E}(\com^{n})$ pour définir $V_{E}$, si l'on suppose que $E$ est compact.
      En effet, si
           $$
               L_{E}^{+}(\com^{n}):=\braces{ u \in PSH(\com^{n}): u\vert_{E}\leqslant 0,\; u\pgoe0,\;
               \exists c_{u} \in \ree , \quad u \leqslant
                c_{u} + \log ( 1 + \vert z \vert) },
           $$
      on a même
                $$
                    V_{E}(z):= \sup \braces{ u(z): u\in L_{E}^{+}(\com^{n})\cap\mathcal{C}^{\infty}(\com^{n})}, \quad
                    \forall z \in \com^{n}.
                $$

       On appellera la régularisée supérieure de la fonction $V_{E}$ la plus petite fonction semi-continue
       supérieurement majorant $V_{E}$. On la note $V_{E}^{\star}$ et on la définit comme ci-dessous:
                $$
                        V_{E}^{\star}(z):=\limsup_{\zeta\rightarrow
                        z,\;\zeta\in\com^{n}}V_{E}(\zeta).
                $$

       Si $E$ est non pluripolaire (\textit{i.e.} $\forall u\in PSH(\com^{n}),\;E\not\subset u^{<-1>}(-\infty)$)
       la fonction $V_{E}^{\star} \in PSH(\com^{n})$. Nous avons le théorème approximation de
       $V_{E}$
       dû à Siciak.
                    \begin{equation}
                             V_{E}= \log (\Phi_{E})
                    \end{equation}
 où
$$ \Phi_{E}(z):=\sup \braces{ \vert p(z) \vert^{\frac{1}{\deg(p)}}:p \in \polyn{\com}{z}{n}, \; \norm{p}_{E}\leqslant 1, \;
 \deg(p) \geqslant 1 }, \quad  z \in \com^{n}.
$$ $\norm{\cdot}_{E}$ est la norme uniforme sur le compact $E$. A partir de l' égalité (2), nous déduisons trivialement
l'inégalité de Bernstein-Walsh :
\begin{equation}
\modul{p(z)} \leqslant \norm{p}_{E}  \, e^{(\deg(p) \,V_{E}(z))}, \quad \forall p \in
\polyn{\com}{z}{n} , \; \forall z \in \com^{n}.
\end{equation}

On trouvera plus de précisions sur les démonstrations des propriétés de la fonction de
Green avec pôle à l'infini et du théorème de Siciak dans son papier fondamental
\cite{siciak} ou dans l'ouvrage \cite{klimek}.

Si $E$ est un compact de $\com^n$, non pluripolaire, nous dirons que $E$ à la propriété
d'$HCP$ ($H$ölder $C$ontinuity $P$rincipe), s'il existe $\delta_0,\;C,\;\alpha>0$ des
constantes ne dépendant que de $E$ telles que:
$$
          \forall \delta\in\bracks{0,\delta_0},\;\forall z\in\com^{n},\;d(z,E)\ppoe\delta\Longrightarrow V_E(z)\ppoe
          C\delta^{\alpha}\textrm{ ($d$ métrique Euclidienne)}.
$$

 Si un compact à la propriéte d'$HCP$, la fonction $V_E$ est continue; ceci implique
que le compact $E$ est $L$-régulier (c.f.\cite{siciak}).

Nous dirons que $E$ est localement $HCP$, si pour tout $a\in E$, il existe $\eps>0$, tel
que $E\cap B(a,\eps)$ vérifie la propriété d'$HCP$. Pour plus de commodité grammaticale
et syntaxique, nous écrirons simplement que $E$ est $HCP$ ou $E$ est \textit{localement}
$HCP$. Parfois, par abus de language, nous préciserons la valeur de la constante $\alpha$
en parlant d'$HCP$ d'exposant ou d'ordre $\alpha$, s'il n'y a pas de confusions possible.

Il est évident, qu'avec l'inégalité de Benstein-Walsh (3) et les inégalités de Cauchy,
qu'un compact de $E\subset\com^{n}$  $HCP$ admette des Inégalités de Markov.

La notion de d'$HCP$ a été explorée  par Plesniak et Pawlucki dans l'article
\cite{pawlucki2} de manière très approfondie.

\section{Continuité de Hölder de la fonction extrémale dans le sous-ensemble algébrique $\tilde{A}$}

Dans cette partie, nous allons affaiblir la notion d'$HCP$ pour les courbes algébriques
de $\ree^{n}$ $(n\pgoe 2)$ en la restreignant au sous-ensemble algébrique complexifié.
Pour ce faire, nous introduirons la métrique des géodésiques dans la courbe complexifiée
et nous verrons, modulo une certaine compatibilité pour cette métrique, que cette notion
d'$HCP$ impliquera encore l'existence d'inégalités de Markov tangentielles.

En définitive, il ne sera pas nécessaire de montrer la propriété d'$HCP$ dans $\com^{n}$,
mais seulement dans une partie de dimension complexe 1.

\subsection{Prolongement de la paramétrisation de Puiseux pour les courbes algébriques complexes}

Pour obtenir le prolongement de la paramétrisation de Puiseux, il est nécessaire de
construire une projection globale $\pi:\tilde{A}\rightarrow \com_{1}$ propre.
\begin{pro}
      Soient $G=\prim{G}\times\prims{G}$, où $\prim{G}$ et $\prims{G}$ sont
      deux sous-ensembles ouverts respectifs de $\com^{p}$ et $\com^{m}$
      $(m+p=n)$, et $\pi:(\prim{z},\prims{z})\mapsto \prim{z}$. Soit
      $\tilde{A}$ un sous-ensemble analytique de $G$ tel que $\pi : \tilde{A} \longrightarrow
      \prim{G}$ soit une application propre. Alors
      $\prim{\tilde{A}}=\pi(\tilde{A})$ est un sous-ensemble analytique de
      $\prim{G}$, et le nombre de pré-images $\pi^{-1}(\braces{\prim{z}}),\;(\prim{z}\in\prim{G})$ est
      localement fini dans $\prim{G}$. Si, de plus, $G=\com^{n}$,
      $\prim{G}=\com^{p}$ et $\tilde{A}$ est un sous-ensemble algébrique de
      $\com^{n}$, alors $\prim{\tilde{A}}=\pi(\tilde{A})$ est aussi un
      sous-ensemble algébrique de $\com^{p}$.
\end{pro}
On trouvera la démonstration de la proposition 2 dans le livre de Chirka (\cite{chirka},
\S 3.2, page 29).

\begin{lem}
Soit $A$ un sous-ensemble algébrique de $\ree^{n}$ de dimension 1 et $\tilde{A}$ son
complexifié dans $\com^{n}$
  tel que $A=\tilde{A} \cap \ree^{n}$. Alors il existe l une transformation unitaire de $\com^{n}$
   telle que $l(\ree^{n}) = \ree^{n}$ et
$\pi : l(\tilde{A}) \longrightarrow \com_{1}$ soit propre.
\end{lem}

\textbf{Preuve:} On injecte $\com^{n}$ dans $\proj_{n}(\com)$, $H_{0}$ sera l'hyperplan à
l'infini identifié $\proj_{n-1}(\com)$ dans $\proj_{n}(\com)$ de sorte que
$\proj_{n}(\com)=\proj_{n-1}(\com)\cup\com^{n}$. La courbe  $A$ est algébrique, donc il
existe des polynômes $p_{1},\ldots,p_{s}$ dans $\polyn{\ree}{x}{n}$, $s\in\nat^{\star}$,
tels que:
  $$
    A=\braces{x\in\ree^{n}:p_{1}(x)=\cdots =p_{s}(x)=0}.
  $$
Si $d_{j}$ est le degré du polynôme $p_{j}$, pour tout $j\in\braces{1,\ldots,s}$, nous
avons alors la décomposition suivante pour tous les polynômes $p_{j}$:
   $$
     p_{j}(x)=h_{j}(x)+q_{j}(x),\quad \forall j \in \braces{1,\ldots,s},
   $$
où les $h_{j}$ sont des polynômes homogènes dans $\polyn{\ree}{x}{n}$ de degré $d_{j}$ et
les $q_{j}$ sont des polynômes dans $\polyn{\ree}{x}{n}$ tels que $\deg(q_{j})<d_{j}$.
Considérons le sous-ensemble algébrique projectif $V$, défini comme ci-dessous:
  $$
    V:=\braces{[z]\in\proj_{n-1}(\com):h_{1}(z)=\cdots=h_{s}(z)=0}
  $$
  avec
  $$
    V\subset H_{0},\quad H_{0}=\proj_{n-1}(\com).
  $$
Le sous-ensemble algébrique projectif $V$ est bien défini puisque les $h_{j}$ sont des
polynômes homogènes et $V$ est propre, car les $h_{j}$ ne sont pas tous nuls. Maintenant
identifions $Gr_{\com}(n,1)$ et $\proj_{n-1}(\com)$. Soit
  $$
    \mathcal{L}:=\braces{\tilde{L}\in Gr_{\com}(n,1):\tilde{L}\subset
    \tilde{A}},
  $$
il est évident de voir que $\mathcal{L}$ est un fermé d'intérieur vide dans
$Gr_{\com}(n,1)$. Posons donc
  $$
     \mathcal{M}:=Gr_{\com}(n,1)\setminus \mathcal{L},
  $$
il s'en suit que $\mathcal{M}$ est un ouvert partout dense dans $Gr_{\com}(n,1)$.
Définissons l'ensemble suivant:
  $$
     \mathcal{B}:=\braces{\tilde{L}\in Gr_{\com}(n,1):\exists v_{n}\in \ree^{n},\;\norm{v_{n}}_{2}=1,\; \tilde{L}=\com v_{n},\;
     [v_{n}]
      \not\in V}.
  $$
De toute évidence $\mathcal{B}\not=\emptyset$, car $\deg(h_{j})=\deg(p_{j})$ implique que
$V$ est au plus une hypersurface. Supposons maintenant que
$\mathcal{B}\not\subset\mathcal{M}$. Donc il existe $\tilde{L}\in \mathcal{B}$ tel que
$\tilde{L}\not \in \mathcal{M}$. Il s'en suit que $\tilde{L}\subset \tilde{A}$.
Choisissons $v_{n}$ dans $\ree^{n}$ tel que $\tilde{L}=\com v_{n}$, car $\tilde{L}$ est
dans $\mathcal{B}$. Posons $L=\ree v_{n}$. Nous avons:
   $$
     L\subset\tilde{L}\subset\tilde{A}\Longrightarrow L\subset A.
   $$
Donc
  $$
     \forall \lambda \in \ree,\quad p_{1}(\lambda v_{n})=\cdots=p_{s}(\lambda
     v_{n })=0
  $$
Soit encore
  $$
    \forall \lambda >0,\quad \frac{1}{\lambda^{d_{1}}}p_{1}(\lambda v_{n})=\cdots=\frac{1}{\lambda^{d_{s}}}p_{s}(\lambda
     v_{n})=0.
  $$
Si $\lambda \rightarrow \infty$, cela nous donne:
  $$
    h_{1}(v_{n})=\cdots=h_{s}(v_{n})=0,
  $$
donc $[v_{n}]$ est dans $V$, ce qui est une contradiction. En conclusion
$\mathcal{B}\subset\mathcal{M}$. Soient, maintenant, les polynômes homogénéisés des
$p_{j}$ définis comme ci-dessous:
   $$
     p_{j}^{\star}(z_{0},\ldots,z_{n}):=z_{0}^{d_{j}}p_{j}\parentheses{\frac{z_{1}}{z_{0}},\ldots,\frac{z_{n}}{z_{0}}},\quad
      j \in\braces{1,\ldots, s}.
   $$
 Les $p_{j}^{\star}$ sont des polynômes homogènes de $\com\lbrack z_{0},\cdots,
z_{n}\rbrack$. Considérons la variété projective de $\proj_{n}(\com)$
  $$
    \Tilde{\Tilde{A}}:=\braces{[z]\in\proj_{n}(\com): p_{1}^{\star}(z)=\cdots=\
    p_{s}^{\star}(z)=0}.
  $$
Comme pour tout $j\in \braces{1,\ldots,s}$, on a $p_{j}=p_{j}^{\star}\vert_{\com^{n}}$.
Il s'en suit que $\tilde{A}=\com^{n}\cap\Tilde{\Tilde{A}}$. Étant donné que
$\Tilde{\Tilde{A}}$ est fermé, nous en déduisons successivement les inclusions suivantes:
  $$
    \overline{\Tilde{A}}\subset\Tilde{\Tilde{A}} \;\textrm{ et }\;
    \overline{\Tilde{A}}\setminus\tilde{A}\subset\Tilde{\Tilde{A}}.
  $$
  On a aussi $\Tilde{\Tilde{A}}\cap H_{0}=V$, par conséquent $\overline{\tilde{A}}\cap H_{0}\subset V$, et $\tilde{A}\cap
  H_{0}=\emptyset$. Avec l'inclusion ci-dessous:
    $$
      \tilde{A}\subset\overline{\tilde{A}}\subset\Tilde{\Tilde{A}}
    $$
et
    $$
      \tilde{A}\subset\overline{\tilde{A}}\cap
      \com^{n}\subset\Tilde{\Tilde{A}}\cap\com^{n}=\tilde{A}.
    $$
Donc
    $$
      \tilde{A}=\overline{\tilde{A}}\cap
      \com^{n}=\Tilde{\Tilde{A}}\cap\com^{n}.
    $$
Montrons que $\overline{\tilde{A}}\setminus\tilde{A}\subset V$. On a
   $$
     \overline{\tilde{A}}=\overline{\tilde{A}}\cap\parentheses{H_{0}\cup\com^{n}}=\parentheses{\overline{\tilde{A}}\cap
     H_{0}}\cup\parentheses{\overline{\tilde{A}}\cap\com^{n}}=\parentheses{\overline{\tilde{A}}\cap
     H_{0}}\cup\tilde{A}.
   $$
De plus
$\overline{\tilde{A}}=\parentheses{\overline{\tilde{A}}\setminus\tilde{A}}\cup\tilde{A}$.
Donc on a l'égalité ci-dessous:
   $$
     \parentheses{\overline{\tilde{A}}\setminus\tilde{A}}\cup\tilde{A}=\parentheses{\overline{\tilde{A}}\cap
     H_{0}}\cup\tilde{A}\quad (\Diamond_{1})
   $$
Choisissons $[z]\in\overline{\tilde{A}}\setminus\tilde{A}$, on a $[z]\not\in\tilde{A}$.
Avec l'égalité $(\Diamond_{1})$, on en déduit que $[z]\in\overline{\tilde{A}}\cap
     H_{0}$, soit $\overline{\tilde{A}}\setminus\tilde{A}\subset\overline{\tilde{A}}\cap
     H_{0}$. L'inclusion voulue est démontrée, sachant que nous avons $\overline{\tilde{A}}\cap H_{0}\subset
     V$, d'où:
     $$
     \overline{\tilde{A}}\setminus\tilde{A}\subset V\subset
     H_{0}.
     $$
      Choisissons maintenant $\tilde{L}$ dans
$\mathcal{B}$. Nous avons $\tilde{L}\not\subset\tilde{A}$ et $L\not\subset A$. Il s'en
suit que $\tilde{L}\cap\tilde{A}$ est une sous-ensemble algébrique propre de $\tilde{L}$.
Comme $\tilde{L}$ est de dimension complexe 1 et $\tilde{A}$ algébrique,
$\tilde{L}\cap\tilde{A}$ est fini, donc $L\cap A$ aussi. Construisons maintenant la
transformation unitaire $l$. Choisissons $v_{n}$ tel que $\norm{v_{n}}_{2}=1$ et $\com
v_{n}=\tilde{L}$. Rappelons que $\com^{n}$ est muni du produit scalaire ci-dessous:
  $$
    \forall z,\zeta\in \com^{n},\quad(z\vert
    \zeta):=\sum_{j=1}^{j=n}z_{j}\overline{\zeta}_{j},
  $$
avec $z=\uplet{z}{n},\;\zeta=\uplet{\zeta}{n}$.
 Construisons par le procédé d'orthonormalisation de Gram-Schmidt une base
$\ree$-orthogonale dans $L^{\bot_{\ree}}$. Nommons la $(v_{1},\ldots,v_{n-1})$. La base
$(v_{1},\ldots,v_{n-1})$ est $\com$-libre. Nous avons donc $\dim_{\com}
Vect_{\com}(v_{1},\ldots,v_{n-1})=n-1$. Il s'en suit que :
  $$
    Vect_{\com}(v_{1},\ldots,v_{n-1})=\tilde{L}^{\bot_{\com}}\text{ tel que }\tilde{L}\oplus\tilde{L}^{\bot_{\com}}=\com^{n}.
  $$
  On a
  $$
    \forall a\in \tilde{A},\quad\parentheses{\overline{\tilde{A}}\setminus\tilde{A}}\cap\overline{\tilde{L}}
    \subset V \cap\overline{\parentheses{a+\tilde{L}}}=\emptyset,
  $$
  car $V$ est un sous-ensemble de $H_{0}$ et $[v_{n}]$ n'appartient pas à $V$. On en déduit que la projection ci-dessous:
  $$
     \begin{array}{ccccc}
     \pi_{\tilde{L}} &:& \tilde{A} & \longrightarrow & \tilde{L}^{\bot_{\com}}\quad (\Diamond_{2})\\
     \end{array}
  $$
  est propre. D'après la construction de $\tilde{L}$ et $\tilde{L}^{\bot}$,
  il existe donc une unique transformation unitaire $l_{1}$ de $\com^{n}$ telle
que
  $$
    (l_{1}(v_{j})\vert e_k)=\delta_{j,k},\quad \forall j,\,k \in
    \braces{1,\ldots,n},
  $$
  avec $\delta_{i,j}$ symbole de Kronecker, et on a
  $$
    l_{1}(\ree^{n})=\ree^{n}.
  $$
 On déduit d'après ce qui précède, c'est-à-dire $(\Diamond_{2})$, que la  projection :
  $$
   \begin{array}{ccccc}
   \pi^{1} & : & l_{1}(\tilde{A}) & \longrightarrow & \com_{n-1}\\
       &   & (\prim{z},z_{n})&\longmapsto & \prim{z}\\
   \end{array}
  $$
est propre. Notons maintenant $\tilde{A}_{1}=\tilde{A}$. D'après la proposition 3,
$\pi^{1}(l(\tilde{A}_{1}))$ est algébrique dans $\com_{n-1}$. Construisons par récurrence
descendante la suite $(\pi^{i},l_{i},\tilde{A}_{i})$, pour tout $i\in \braces{2,\ldots,
n-1}$, telle que:
   $$
     \pi^{i} :l_{i}(\tilde{A}_{i}) \longrightarrow
     \com_{n-i},\quad\tilde{A}_{i}=\pi^{i-1}\parentheses{l_{i-1}(\tilde{A}_{i-1})},\quad
     \forall i \in\braces{2,\ldots, n-1},
   $$
où $\pi^{i}$ est une projection propre telle que $\pi^{i}(\ree^{n-i})=\ree^{n-i-1}$,
$\tilde{A}_{i}$ est un sous-ensemble algébrique et $l_{i}$ est dans $O(\com,n-i)$,
l'ensemble des transformations unitaires de $\com^{n-i}$ et
$l_{i}(\ree^{n-i})=\ree^{n-i}$. Nous avons les inclusions suivantes:

   $$
      \begin{array}{ccccc}
       j_{i} & : & O(\com,i)& \longrightarrow & O(\com,i+1)\\
             &   & l_{n-i+1} & \longmapsto & l_{n-i+1}\oplus id_{\com_{(i+1)}}
       \end{array}
       , \quad \forall i\in \braces{2,\ldots,n-1},
   $$
   où $\com_{(i)}$ est le sous-espace $\braces{0_{\com^{i-1}}}\times\com$ de $\com^{n}$.
 Il s'en suit que si l'on pose $\tilde{l}=(l_{n-1}\oplus
id_{\com_{(3)}})\circ\cdots\circ(l_{2}\oplus id_{\com_{(n)}})\circ l_{1}$, $\tilde{l}$
 se trouve dans $O(\com,n)$, telle que $\tilde{l}(\ree^{n})=\ree^{n}$. Si
$\pi=\pi^{n}\circ\cdots\circ\pi^{1}$, où $\pi^{n}:
l_{n-1}(\tilde{A}_{n-1})\longrightarrow \com_{1}$ est la projection propre choisie comme
$\pi^{1}$ pour le sous-ensemble algébrique $\tilde{A}_{n-1}$. On a donc :
  $$
   \begin{array}{ccccc}
   \pi & : & l(\tilde{A}) & \longrightarrow & \com_{1}\\
       &   & \uplet{z}{n}&\longmapsto & z_{1}\\
   \end{array}
  $$
est propre car toutes les projections $\pi^{i}$ le sont, donc $\pi$ et $l$
 vérifient les conditions requises du lemme.
 $\blacksquare$

\begin{pro}
                Supposons $A$ une courbe algébrique réelle,
                localement irréductible, telle que
                 $0 \in \tilde{A}_{sing}$ et $\varphi$ la
                paramétrisation de Puiseux de la proposition 2.
                Modulo un changement de
                coordonnées unitaire dans $\com^{n}$, laissant $A$ dans $\ree^{n}$,
                 il existe $\Psi$ une application
                holomorphe définie dans $\prim{\Omega}$ ouvert de
                $\com$ partout dense, avec $\overline{D(0,1)}\subset
                \prim{\Omega}$, telle que $\Psi$ prolonge $\psi$ à
                $\prim{\Omega}$ et pour tout $z$ dans
                $\prim{\Omega}$, $\modul{\Psi(z)}\leqslant
                C(1+\modul{z}^{sk})$, où $C,s>0$ sont des
                constantes ne dépendant que de $\tilde{A}$.
       \end{pro}

 \textbf{Preuve:} On peut donc choisir, d'après le lemme 2 du 4.1, un changement de coordonnées unitaire et une projection
 $\pi : \tilde{A}
\longrightarrow \com_{1}, \quad \pi \uplet{\xi}{n}= \xi_{1},$ propre telle que
$\pi(\ree^{n})=\ree$. Comme $\tilde{A}$ est supposée algébrique, d'après la proposition 1
du 3.2, il existe deux constantes réelles strictement positives $C,s$ telles que:
\begin{equation}
( \modul{\xi_{2}}^{2}+ \cdots +  \modul{\xi_{n}}^{2} )^{\frac{1}{2}}\leqslant
C(1+\modul{\xi_{1}}^{s}).
 \end{equation}
La projection $\pi$ est propre, donc il existe $\sigma_{1} \subset \com_{1}$ un
sous-ensemble fini ($0\in\sigma_1$ car $0\in \tilde{A}_{sing}$), sachant que  $\tilde{A}$
est de dimension complexe 1 et algébrique, tel que,
 $$ \pi :
\tilde{A} \setminus \pi^{<-1>}(\sigma_{1}) \longrightarrow \com \setminus \sigma_{1}
 $$
soit un $r$-revêtement holomorphe ($r \in \nat^{\star}$), avec
$$
 \forall z_{1} \in \com
\setminus \sigma_{1}, \quad \textrm{card}\;\parentheses{\pi^{<-1>}(z_{1}) \cap
\tilde{A}}=r.
$$
 Donc pour tout $z_{1} \in \com_{1} \setminus \sigma_{1}$, il existe de $V_{z_{1},j}$, ($ \forall j  \in \braces{1,\ldots,r}$)
 des ouverts de $\tilde{A} \setminus \pi^{<-1>}(\sigma_{1})$ et $W_{z_{1}} \in \mathcal{V}_{\com_{1}}(z_{1})$
(voisinage ouvert de $z_{1})$ avec $\alpha_{z_{1},j} \in H(W_{z_{1}})  \quad ( \forall j
\in \braces{1,\ldots,r}).$
$$
\begin{array}{cccc}
\alpha_{z_{1},j} : & W_{z_{1}} &  \longrightarrow    & V_{z_{1},j}
\\
                           & \xi              &  \longmapsto         & (\xi,\alpha_{z_{1},j }(\xi )). \\
\end{array}
$$ Posons $\Omega = \com_{1} \setminus \Delta$, où $\Delta$ est un ensemble
fini de demi-droites de $\com_{1}$ ayant pour origine tous les points de $\sigma_{1}$,
contenant en particulier $\cobracks{0,+\infty}$, rendant ainsi $\Omega$ simplement
connexe\footnote{L'auteur remercie Julien Duval pour cette idée.}. En appliquant le
théorème de la monodromie, il existe $\alpha_{1}, \ldots , \alpha_{r} \in H(\Omega)$ tels
que : $$ \forall z_{1} \in \Omega, \ \exists W_{z_{1}} \in\mathcal{V}_{\com_{1}}(z_{1}),\
\forall j \in \braces{1,\ldots,r}, \quad \alpha_{j} \vert_{W_{z_{1}}}= \alpha_{z_{1},j}.
$$
 Modulo un changement d'indices choisissons $U \in \mathcal{V}_{\com^{n}}(0)$ où
$U=U^{\prime} \times U^{\prime \prime} \subset \com \times \com^{n-1}$
 tel que:
$$
\begin{array}{cccc}
\pi\vert_{ U}: & \tilde{A} \cap U & \longrightarrow & U^{\prime}
\\
                     & \uplet{\xi}{n} & \longmapsto & \xi_{1} \\
\end{array}
$$ soit une application propre, $\pi^{<-1>}(0) \cap \tilde{A} \cap U= \lbrace 0
\rbrace$ et $ \tilde{A} \cap U$ soit irréductible, rappelons que la courbe algébrique
réelle $A$ est localement irréductible dans les hypothèses de notre lemme.
 D'après le théorème de structure locale des sous-ensembles analytiques, il existe $\sigma_{2} \subset \com_{1}$
 ($\sigma_2\subset\sigma_1$) fini tel que $\pi : (\tilde{A} \cap U) \setminus \pi^{<-1>}(\sigma_{2}) \longrightarrow
\com_{1} \setminus \sigma_{2}$
 soit un $k$-revêtement holomorphe ramifié ($k \in \nat^{\star}$).
Construisons la paramétrisation de Puiseux associée à la singularité $0$ de $\tilde{A}$ à
partir de la projection $\pi\vert_{U}$ (c.f. la proposition 2 du 3.4), soit $\varphi(z) =
(c z^{k},\psi (z))$ cette paramétrisation avec $k=\mu_0(\tilde{A})$. Considérons les
secteurs suivants, si $\overline{D^{\prime}(0,1)}=\overline{D(0,1)}\setminus \lbrace 0
\rbrace$,
    $$
             \mathcal{S}_{j}:= \left \lbrace
             z \in \prim{D}(0,1): \frac{2  \pi j}{k}< \arg (z)\ppoe \frac{2  \pi (j+1)}{k} \right \rbrace, \quad
             (\forall j \in \braces{0,\ldots,k-1}).
    $$
 Définissons les changements de coordonnées locales suivantes :
    $$
             \begin{array}{cccc}
                    \theta_{j}: & \overline{D^{\prime}(0,1)} & \longrightarrow & \mathcal{S}_{j} \\
                                 & \xi                    & \longmapsto     & (\frac{\xi}{c})^{\frac{1}{k}} \\
             \end{array}.
     $$
 On a évidemment $\theta_{j}
 \in H(\overline{\prim{D}(0,1)})$. Posons $\prim{\mathcal{S}_{j}} =\mathcal{S}_{j} \setminus \Delta_{j} ,\;
 \Delta_{j}:=\theta_{j} (\Delta)$, et $\Omega_{j}:=\theta_{j}
(\Omega)$, on voit que $\Omega_j\subset\prim{\mathcal{S}}_j$. D'après le théorème de
l'image ouverte $\Omega_{j}$ est ouverte. Il est clair que $\theta_{j}$  paramétrise
$\prim{\mathcal{S}}_{j}$ sur $\overline{D(0,1)} \setminus
\Delta\subset\overline{\prim{D}(0,1)}$. Donc $\varphi \circ \theta_{j}$ paramétrise un
morceau d' une feuille au dessus de $\overline{D(0,1)} \setminus \Delta\subset\com_1$. On
en déduit l'existence d'un $\alpha_{i_{j}}$, tel que $\psi \circ \theta_{j}=
\alpha_{i_{j}}$ sur $\Omega$. Donc la fonction $\psi$ se prolonge sur l'ouvert
$\Omega_{j}$ et de l'inégalité (4) on déduit, pour tout $z \in \Omega_{j}, \quad
\modul{\psi(z)} \leqslant c_{j}(1+\modul{z})^{sk} \quad (\forall j \in
\braces{0,\ldots,k-1})$. On déduit de même que $\psi$ se prolonge sur chaque
$\Omega_{j}$. Comme les prolongées coïncident sur $\overline{\prim{D}(0,1)} \setminus
\prim{\Delta}$, alors d'après le théorème du prolongement analytique, $\psi$ se prolonge
sur l'ouvert $\prim{\Omega}$ défini par :
    $$
           \prim{\Omega}= \bigcup_{1\leqslant j \leqslant k-1 } \Omega_{j}.
    $$
 Soit $\Psi$ sa prolongée, elle vérifie donc:
    $$
           \forall z \in \prim{\Omega}, \quad \modul{\Psi (z)} \leqslant c(1+
           \modul{z}^{sk}),
    $$
 avec $\overline{\prim{\Omega}}=\com_{1}.\;\blacksquare$

Dans le corollaire 3, nous démontrons, sans le dire, une version globale du théorème des
fonctions implicites.

Précisons dans le corollaire ci-dessous ce qui diffère lorsque l'on choisit le point
$x_{0}$ dans $A_{reg}$.
         \begin{cor}
               Soit $A$ une courbe algébrique réelle de $\ree^{n}$
               localement irréductible. Soit $x_{0}$ dans $\tilde{A}_{reg}$ fixé. Modulo un changement de coordonnées unitaire
                laissant $A$ dans $\ree^{n}$, il existe
               une paramétrisation holomorphe $\varphi(z)=(cz,\psi(z))$, telle que $\varphi:\overline{D(0,1)}\rightarrow \tilde{A}\cap
               U$, où $U=\varphi(\overline{D(0,1)})$ et $\varphi([-1,1])=A\cap U$. Cette paramétrisation se
               prolonge dans $\Omega$ un ouvert de $\com$ partout dense, telle que
               $\modul{\psi(z)}\ppoe c(1+\modul{z}^{s})$, pour tout $z$ dans $\com$. Où $c$ et $s$ sont des constantes
               dépendant uniquement de la courbe $A$.
         \end{cor}
    Étant donné la forte similitude de cette démonstration avec celle du corollaire 3,
     nous nous permettrons d'être un peu plus expéditif, quant à la démonstration de celui-ci.

     \textbf{Preuve:} D'après le lemme 2 du 4.1, il existe un changement de coordonnées unitaire laissant
     $A$ dans $\ree^{n}$, tel que la projection,
         $$
              \begin{array}{ccccc}
              \pi & : & \tilde{A}    & \longrightarrow & \com_{1} \\
                  &   & \uplet{z}{n} & \longmapsto     & z_{1}\,,\\
              \end{array}
         $$
     soit propre. Sans perdre de généralité, on peut supposer que $x_{0}=0$ et $0\in A_{reg}$. Il
     existe $\sigma$ un sous-ensemble, formé de points isolés, inclus dans $\com$, tel que
         $$
              \pi:\tilde{A}\setminus\pi^{<-1>}(\sigma)\longrightarrow\com\setminus\sigma,
         $$
     soit un $r$-revêtement ramifié. En opérant de la même manière
     qu'à la proposition 4 du 4.1, il existe un ouvert $\Omega$ de
     $\com$ partout dense et simplement connexe et des fonctions $\alpha_{j}$
     ($j\in\braces{1,\ldots,r}$) holomorphes dans
     $\Omega$ telle que $\psi=\alpha_{j}\vert_{D(0,1)}$. Étant donné que la projection est
     propre et la courbe $\tilde{A}$ est algébrique, on a, d'après la proposition 1 du 3.2,
     l'estimation suivante:
         $$
             \modul{\psi(z)}\ppoe c(1+\modul{z}^{s}),\quad\forall z\in \Omega,
         $$
         d'où le corollaire souhaité. $\blacksquare$
\subsection{Estimation de la fonction extrémale}
Nous allons estimer la fonction de Green avec pôle à l'infini sur un voisinage de la
courbe algébrique réelle. Nous remarquerons que la proposition 5, ci-après, nous donne
une autre démonstration d'une partie du théorème de Sadullaev \cite{sadullaev}.

Dans la proposition suivante les hypothèses sur $A$ sont les mêmes que dans la
proposition 4 du 4.1 et la paramétrisation de Puiseux $\varphi$ est celle obtenue dans la
proposition 2 du 3.4.

\begin{pro}
Modulo un changement de coordonnées unitaire laissant $A$ dans $\ree^{n}$, il existe des
constantes réelles $c_{1},c_{2},\eps_{0}$ et $\rho_{1}$ strictement positives dépendant
uniquement de $\varphi$ telles que: $\forall \eps \in \lbrack 0, \eps_{0}\lbrack, \;
\forall z \in D(0,\rho)$,

\begin{equation}
\begin{array}{c}
 \quad V_{A
\cap B(0,c_{1} \eps^{k}) } \circ \varphi(z) \leqslant c_{2} V_{\eps
\mathcal{R}^{<r_A^+>}}(z),\\
                      \\
\parentheses{\textrm{resp. }V_{A
\cap B(0,c_{1} \eps^{k}) } \circ \varphi(z)
                              \leqslant c_{2} V_{\eps
                          \mathcal{R}^{<r_A^->}}(z),}\\
\end{array}
\end{equation}
avec
   $$
     \mathcal{R}^{<r_A^+>}=\bigcup_{1 \ppoe j \ppoe r_{A}^+} \left\lbrack 0, e^{\frac{2ik_{j}\pi}{k}}
\right\rbrack,\;\parentheses{\textrm{resp. }
     \mathcal{R}^{<r_A^->}=\bigcup_{1 \ppoe j \ppoe r_{A}^-} \left\lbrack 0, e^{\frac{2ik_{j}'\pi}{k}}
\right\rbrack,}
   $$
       où $k=\mu_0(\tilde{A})$, $r_{A}^+$ (resp. $r_{A}^-$) est le nombre de composantes
       connexes de $A\cap B(0,c_{1}\eps^{k})\cap\braces{\uplet{z}{n}\in\com^n:\Re
       e(z_1)>0}$  (resp. le nombre de composantes connexes de $A\cap B(0,c_{1}\eps^{k})\cap
       \braces{\uplet{z}{n}\in\com^n:\Re
       e(z_1)<0}$ ). L'ensemble $\mathcal{R}^{<r_A^+>}$ et $0\ppoe k_{1}<\cdots<k_{r_{A}^+}\ppoe k$ (resp.
           $\mathcal{R}^{<r_A^->}$ et $0\ppoe k_{1}'<\cdots<k_{r_{A}^{-}}'\ppoe k$) sont
           construits comme dans la proposition 2 du 3.4.

\end{pro}

\textbf{Preuve:} Montrons l'estimation (5) pour $\mathcal{R}^{<r_A^+>}$ seulement, la
démonstration sera la même pour $\mathcal{R}^{<r_A^->}$. Posons
$\varphi(z)=(cz^k,\psi_2(z)\ldots,\psi_n(z)),\;$ $\forall z \in\overline{D(0,1)}$ (cf. la
proposition 2 du 3.4); soit $\eps_{0}\in\obracks{0,1}$, de sorte que
$\varphi(\overline{D(0,\eps)})\subset B_{\com^{n}}(0,c_{1} \eps^{k})$, pour tout $\eps$
dans $\obracks{0,\eps_{0}}$. Remarquons que nous avons les inclusions suivantes:
$\forall\eps\in\obracks{0,\eps_{0}}$,
     \begin{equation}
                 \varphi (\mathcal{R}^{<r_A^+>}) \subset A , \quad \varphi ( \eps \mathcal{R}^{<r_A^+>})
                 \subset A \cap \varphi (\overline{D(0, \eps)})  \subset A \cap B_{\com^{n}} (0,c_{1}
                 \eps^{k}).
     \end{equation}
Celles-ci entraînent les inclusions ci-après:
      $$
          L_{A\cap B(0,c_{1}\eps^{k})}(\com^{n})\subset
          L_{A\cap\varphi(D(0,\eps))}(\com^{n})\subset L_{\varphi(\eps\mathcal{R}^{<r_A^+>})}(\com^{n}).
      $$
 Nous allons directement majorer la fonction $V_{A\cap B(0,c_{1}\eps^{k})}$. D'après la proposition 4 du 4.1,
  $\varphi$ se prolonge en une fonction holomorphe dans un ouvert partout dense $\Omega$ de
  $\com$ telle que, $\max_{2\ppoe j\ppoe n}\modul{\psi_j(z)}\ppoe C\parentheses{1+\modul{z}^{s}}$, pour tout
  $z\in\Omega$.
       Etant donné que $\eps\mathcal{R}^{<r_A^+>}$
est un continu et que $\varphi$ est une fonction holomorphe de $\com$ dans $\com^{n}$, il
s'en suit que $\varphi(\eps\mathcal{R}^{<r_A^+>})$ n'est pas $L$-polaire, par conséquent
$A\cap B(0,c_{1}\eps^{k})$ non plus. Donc, d'après (\cite{klimek}, Corollaire 5.2.2, page
193), il existe une constante réelle strictement positive $C_{\varphi,A,\eps}$ telle que:
     $$
            u\circ\varphi(z)\ppoe C_{\varphi,A,\eps}+s\log(1+\modul{z}),\quad\forall z\in \Omega,\;
            \forall u\in L_{A\cap B(0,c_{1}\eps^{k})}(\com^{n}).
     $$

 Soit $u$ dans $ L_{A\cap B(0,c_{1}\eps^{k})}(\com^{n})$, il s'en suit, de ce qui précède, que $u\circ\varphi(z)\ppoe
 C_{\varphi,A,\eps}+s\log(1+\modul{z})$, pour tout $z$ dans $\Omega$. Soit $\rho$ une
 constante réelle dans $\obracks{\eps_{0},1}$ et
 considérons la fonction auxiliaire $H_{\alpha}$, comme suit:
    $$
               H_{\alpha}(z):=\bracks{s\log\parentheses{\frac{\modul{z}-\rho_{1}}{1-\rho_{1}}}+\alpha}^{+},\quad\forall
               z\in \com,\;\forall \alpha>0,
    $$
où
    $$
           x^+=\max(0,x),\quad\forall x\in\ree.
    $$
 La fonction $H_{\alpha}$ est sous-harmonique dans  $\com$, $\forall \alpha>0$ et vérifie les propriétés ci-dessous:
        \begin{itemize}
              \item[$\bullet$] $H_{\alpha}(z)=0, \quad\forall z \in D(0,\rho)$.
              \item[$\bullet$] $H_{\alpha}(z)\pgoe \alpha, \quad\forall z \in \com\setminus D(0,1)$.
              \item[$\bullet$] $H_{\alpha}(z)\ppoe c_{\alpha}+s\log(1+\modul{z}),\quad\forall z\in\com$, $c_{\alpha}$
              étant une constante réelle strictement positive.
        \end{itemize}
 Maintenant, choisissons $\alpha>0$ suffisamment grand de sorte que $H_{\alpha}(z)\pgoe
 u\circ\varphi(z)$, pour tout $z$ dans $\Omega\setminus D(0,1)$.
 Pour cela minorons $\Delta_{\alpha}(z):=H_{\alpha}(z)-u\circ\varphi(z)$ dans $\Omega$,
       $$
           \Delta_{\alpha}(z)\pgoe s\bracks{\log\parentheses{1-\frac{\rho}{\modul{z}}}-\log\parentheses{1+\frac{1}{\modul{z}}}}
           -s\log(1-\rho)+\alpha-C_{\varphi,A,\eps}\,,\quad\forall
           z\in\Omega\setminus D(0,1),
       $$
       $$
           \Delta_{\alpha}(z)\pgoe
           s\log\parentheses{1+\frac{-1-\rho}{\modul{z}+1}}+\alpha-C_{\varphi,A,\eps}-s\log(1-\rho),
           \quad\forall z\in\Omega\setminus D(0,1).
       $$
 Donc, si on choisit $\alpha>0$ assez grand, on a $\Delta_{\alpha}(z)\pgoe
 0$ pour tout $z$ dans $\Omega\setminus D(0,1)$.
D'après la construction de l'ouvert
 $\Omega$ de la proposition 4 du 4.1 et le choix de $\alpha$, nous avons
    $$
          \forall z\in\partial\Omega,\quad \limsup_{\zeta\rightarrow
          z,\zeta\in\Omega}u\circ\varphi(\zeta)\ppoe H_{\alpha}(z),
    $$
 donc d'après (\cite{klimek}, Corollaire 2.9.14, page 69), la fonction:
    $$
         W_{\alpha}(z)=
              \left \lbrace
                           \begin{array}{lr}
                                      \max(u\circ\varphi(z),H_{\alpha}(z)),  &z\in\Omega\\
                                      H_{\alpha}(z),                         &  z
                                      \in\com\setminus\Omega
                           \end{array}\right.
    $$
 est bien définie et sous-harmonique dans $\com$, $\frac{1}{s}W_{\alpha}
 \in L_{\varphi(\eps\mathcal{R}^{<r_A^+>})}(\com^{n})$ et
 $\restric{W_{\alpha}}{D(0,\rho)}=u\circ\varphi$. Il s'en suit que
     $$
             u\circ\varphi(z)\ppoe W_{\alpha}(z)\ppoe
             sV_{\eps\mathcal{R}^{<k>}}(z),\quad\forall z\in D(0,\rho).
     $$
 Nous concluons donc
     $$
          u\circ\varphi(z)\ppoe sV_{\eps\mathcal{R}^{<r_A^+>}}(z),\quad\forall z\in
          D(0,\rho_{1}),\;\forall u\in  L_{A\cap B(0,c_{1}\eps^{k})}(\com^{n}).\quad\blacksquare
     $$
 \textbf{Remarque}: Pour ne pas alourdir la démonstration de la proposition 5, nous n'avons pas détaillé
 l'hypothèse où la multiplicité $k$ prend la valeur $1$. Rappelons que nous avons supposé que $0$ est un point de
 $A_{sing}$, donc nécessairement $k$ est supérieur ou égal à $2$. Ceci dit, l'obtention de l'estimation (5),
 pour les points non singuliers de $\tilde{A}$,
  ne diffère pas quant aux techniques de démonstration utilisées, il faut seulement substituer
  la proposition 4 au corollaire 3 du 4.1. C'est pourquoi, il ne nous semble pas nécessaire de démontrer le corollaire,
ci-dessous.

    \begin{cor}
               Soit $A$ une courbe algébrique localement irréductible telle que $0\in A_{reg}$.
               Modulo un changement de coordonnées unitaires laissant $A$ dans $\ree^{n}$, ils
               existent des constantes $\eps_{0}$, $\rho$, $c_{1}$ et $c_{2}$ réelles strictement positives
               dépendant de la paramétrisation $\varphi(z)=(cz,\psi(z))$ du corollaire 2, telles
               que

               $\bullet$ si $0\in A_{reg}\setminus\partial A$,
                         \begin{equation}
                                         \forall \eps\in [0,\eps_{0}[, \forall z \in D(0,\rho),
                                         V_{A\cap B(0,c_{1}\eps)}\circ\varphi(z)\ppoe
                                         c_{2}V_{[-\eps,\eps]}(z),
                         \end{equation}

               $\bullet$ si $0\in A_{reg}\cap\partial A$,
                         \begin{equation}
                                         \forall \eps\in [0,\eps_{0}[, \forall z \in D(0,\rho),
                                         V_{A\cap B(0,c_{1}\eps)}\circ\varphi(z)\ppoe
                                         c_{2}V_{[0,\eps]}(z).
                         \end{equation}
    \end{cor}
Le lemme ci-dessous est dû à Bos \cite{bos3}.
\begin{lem}
Soient $b$ un nombre complexe appartenant à $ \eps I$, $\modul{b} \not = \eps \quad (I
=\lbrack -1,1 \rbrack \subset
 \com)$ et $r$ strictement positif
tel que $r < \eps- \modul{b}$. Alors : $$ \sup_{D(b,r)}V_{\eps I} \leqslant c \log (1 +
r), $$ où $c = \max \left \lbrace 1, \frac{2}{dist(b,\eps \partial I)}\right \rbrace.$
\end{lem}

\textbf{Preuve:} Sans perdre de généralité, on peut supposer $b$ strictement positif.
Pour démontrer le lemme ci-dessus il suffit d'estimer $\sup_{D(0,r)}V_{\prim{I}}$ où
$\prim{I}=\lbrack 2 b - \eps, \eps \rbrack$. On a : $$ \sup_{D(0,r)}V_{\prim{I}} = \log
\left ( \frac{ir}{\eps - b} +\sqrt{\left ( \frac{ir}{\eps - b} \right )^{2}-1} \; \right
). $$
Nous rappelons que la fonction de Green avec pôle a l'infini sur le segment $\lbrack -1,1
\rbrack$ dans $\com$ est connue : $$ V_{\lbrack -1,1 \rbrack}(z) = \log^{+} \modul{z +
\sqrt{z^{2}-1}}, \quad \forall z \in \com . $$ $\blacksquare$
\begin{lem}
La fonction $V_{\lbrack -1,1 \rbrack}$ vérifie la propriété de continuité de Hölder avec
un exposant $\frac{1}{2}$, c'est-à-dire :
$$ (\forall \delta \in \lbrack 0,1 \rbrack)(\forall z \in \com , \; dist(\lbrack -1,1 \rbrack,z)
\leqslant \delta)\Longrightarrow
   V_{\lbrack -1,1 \rbrack}(z) \leqslant C\delta^{\frac{1}{2}}\quad (HCP).
$$
Où $C$ est une constante réelle strictement positive.
\end{lem}
Nous ne donnerons la démonstration de ce dernier lemme, qui pourra être trouvée dans
\cite{pawlucki1} ou dans le livre de M.Klimek \cite{klimek}.
\subsection{Métrique des géodésiques et continuité de Hölder dans le sous-ensemble algébrique $\tilde{A}$}

Définissons la métrique des géodésiques dans $\tilde{A}$. Considérons l'ensemble
$\parentheses{\mathfrak{S}_{[0,1]},\preccurlyeq}$ des subdivisions de l'intervalle
$[0,1]$ muni de la relation d'ordre suivante: $\forall \sigma,\;\tau
\in\mathfrak{S}_{[0,1]}, \;\sigma \preccurlyeq\tau$ est équivalent à dire que la
subdivision $\sigma$ est moins fine que $\tau$. Maintenant, si $\gamma$ est une fonction
de $[0,1]$ dans $\tilde{A}$ et $\sigma$ est dans $\mathfrak{S}_{[0,1]}$, nous définirons
$V_{\sigma}(\gamma)$ de la sorte:
  $$
    V_{\sigma}(\gamma):=\sum_{j=0}^{j=p-1}\norm{\gamma(t_{j+1})-\gamma(t_{j})}_{2},\quad
    \sigma=(t_{0},\ldots,t_{p}),\quad\sigma\in\mathfrak{S}_{[0,1]},
  $$
  où $\norm{z}_2:=(\sum_{j=1}^{j=n}\modul{z_j}^2)^{\frac{1}{2}},\quad\forall z\in\com^n$.
  Définissons le sous-ensemble des fonctions à variation bornée ci-dessous:
  $$
    CVB^{\tiny{(\xi_{1},\xi_{2})}}([0,1],\tilde{A}):=\braces{\gamma \in \mathcal{C}^{0}([0,1],\tilde{A})
    :\gamma(0)=\xi_{1},\,\gamma(1)=\xi_{2},\;\sup_{\sigma\in\mathfrak{S}_{[0,1]}}V_{\sigma}(\gamma)<+\infty},
  $$
  où $\mathcal{C}^{0}([0,1],\tilde{A})$ est l'ensemble des fonctions continues de $[0,1]$ dans
  $\tilde{A}$.
  On peut donc définir désormais une métrique dans $\tilde{A}$ de cette manière,
  $$
    d(\xi_{1},\xi_{2}):=\inf
    \braces{V(\gamma):\gamma\in
    CVB^{\tiny{(\xi_{1},\xi_{2})}}([0,1],\tilde{A})},\quad \forall
    \xi_{1},\xi_{2}\in\tilde{A},
  $$
  où
  $$
    V(\gamma):=\sup\braces{V_{\sigma}(\gamma):\sigma\in\mathfrak{S}_{[0,1]}}.
  $$
  Il est clair que $d(\cdot,\cdot)$ définit bien une métrique dans $\tilde{A}$.
  \begin{lem}
            Soit $\tilde{A}$ une courbe algébrique complexe localement irréductible
            . On considère l'espace métrique $(\tilde{A},d)$, où $d(\cdot,\cdot)$
            est la métrique des géodésiques dans $\tilde{A}$.
             Alors, pour tout $\xi_{0}$ dans $\tilde{A}_{reg}$,
            il existe un voisinage ouvert $U$ de $\xi_{0}$ dans $\tilde{A}_{reg}$, une
            paramétrisation de Puiseux $\varphi:D(0,1)\rightarrow U, \;\varphi(0)=\xi_{0}$
            et des constantes réelles
            strictement positives $c_{1}$ et $c_{2}$, ne dépendant que de $\varphi$ tels que,
            \begin{equation}
               c_{1}\modul{\hat{z}_{1}-\hat{z}_{2}}\ppoe
               d(\varphi(\hat{z}_{1}),\varphi(\hat{z}_{2}))\ppoe
               c_{2}\modul{\hat{z}_{1}-\hat{z}_{2}},\quad\forall \hat{z}_{1},\hat{z}_{2}\in D(0,1).
            \end{equation}
  Si de plus $\xi_{0}$ est dans $\tilde{A}_{sing}$, il existe $U$ un voisinage ouvert
  de $\xi_{0}$ dans $\tilde{A}$, $\tilde{A}_{sing}\cap U =\braces{\xi_{0}}$ et $\varphi:D(0,1)\rightarrow U,
  \;\varphi(0)=\xi_{0}$
  la paramétrisation de Puiseux, on a similairement,
            \begin{equation}
               c_{1}\modul{\hat{z}_{1}}^{k}\ppoe
               d(\varphi(0),\varphi(\hat{z}_{1}))\ppoe
               c_{2}\modul{\hat{z}_{1}}^{k},\quad\forall \hat{z}_{1}\in D(0,1),
            \end{equation}
            où $k$ est la multiplicité complexe du point
  singulier $\xi_{0}$ de $\tilde{A}$.
  \end{lem}

\textbf{Preuve:} Commençons par démontrer (9). Soient $\xi_{0}\in \tilde{A}_{reg}$ et
$\varphi:D(0,1)\rightarrow U$ la paramétrisation de Puiseux. Fixons $\hat{z}_{1}$ et
$\hat{z}_{2}$ dans $D(0,1)$. Choisissons $\gamma$ dans
$CVB^{\tiny{(\varphi(\hat{z}_{1}),\varphi(\hat{z}_{2}))}}([0,1],\tilde{A})$ et
$\sigma\in\mathfrak{S}_{[0,1]}$, telle que $\sigma=(t_{0},\ldots,t_{p})$. Sans perdre de
généralité, on peut supposer que $\xi_{0}=0$. Posons $\xi_{k}=\gamma(t_{k}),\;\forall k
\in\braces{0,\ldots,p}$. Comme $\varphi$ est surjective, ils existent
$z_{0},\ldots,z_{p}$ dans $D(0,1)$ tels que $z_{0}=\hat{z}_{1}$, $z_{p}=\hat{z}_{2}$ et
$\xi_{k}=\varphi(z_{k}),\;\forall k \in \braces{0,\ldots,p}$. Nous avons l'égalité
suivante:
  $$
    \xi_{k+1}-\xi_{k}=\varphi(z_{k+1})-\varphi(z_{k})=\int_{0}^{1}\frac{d}{dt}\varphi(tz_{k+1}+(1-t)z_{k})dt.
  $$
D'où
 $$
   \norm{\gamma(t_{k+1})-\gamma(t_{k})}_{2}=\modul{z_{k+1}-z_{k}}\norm{\int_{0}^{1}\prim{\varphi}(tz_{k+1}+(1-t)z_{k})dt}_{2}
 $$
et donc
   $$
            \norm{\gamma(t_{k+1})-\gamma(t_{k})}_{2}=\modul{z_{k+1}-z_{k}}\parentheses{\sum_{j=1}^{j=n}
            \modul{\int_{0}^{1}\prim{\varphi_{j}}(tz_{k+1}+(1-t)z_{k})dt}^{2}}^{\frac{1}{2}},
   $$
   où les $\varphi_j$ sont les fonctions composantes de $\varphi$.
Étant donné que $\xi_{0}$ n'est pas un point singulier, il existe une constante
strictement positive $c_{1}$ dépendant de $\varphi$ telle que:
  $$
   \parentheses{\sum_{j=1}^{j=n}
   \modul{\int_{0}^{1}\prim{\varphi_{j}}(tz_{k+1}+(1-t)z_{k})dt}^{2}}^{\frac{1}{2}}\pgoe
   c_{1}.
  $$
Donc on a
  $$
    V_{\sigma}(\gamma)\pgoe c_{1}\modul{z_{p}-z_{0}}=
    c_{1}\modul{\hat{z}_{2}-\hat{z}_{1}},\quad\forall
    \sigma\in\mathfrak{S}_{[0,1]},\;\forall \gamma\in
    CVB^{\tiny{(\varphi(\hat{z}_{1}),\varphi(\hat{z}_{2}))}}([0,1],\tilde{A}).
  $$
Il en résulte que
  $$
    d(\varphi(\hat{z}_{1}),\varphi(\hat{z}_{2}))\pgoe
    c_{1}\modul{\hat{z}_{1}-\hat{z}_{2}}.
  $$

Pour l'autre inégalité, il suffit de voir, par définition de la métrique
$d(\cdot,\cdot)$, que l'on a
  $$
    d(\varphi(\hat{z}_{1}),\varphi(\hat{z}_{2}))\ppoe
    \int_{0}^{1}\norm{\frac{d}{dt}\varphi(t\hat{z}_{1}+(1-t)\hat{z}_{2})}_{2}dt
    \ppoe c_{2}\modul{\hat{z}_{1}-\hat{z}_{2}},
  $$
 car $\xi_{0}$ est un point régulier de $\tilde{A}$, d'où l'inégalité (9). Montrons maintenant
 l'estimation (10). Supposons que $\xi_{0}$ est dans $\tilde{A}_{sing}$. On peut supposer sans
 perdre de généralité que $\xi_{0}=0$. Soit $\varphi$ la paramétrisation de Puiseux telle que
 $\varphi:\overline{D(0,1)}\rightarrow U\subset\tilde{A}$ et $z\mapsto(cz^{k},\psi_{2}(z),\ldots,\psi_{n}(z))$,
  où $U$ est un ouvert de $\tilde{A}$, tel que
 $U\cap\tilde{A}_{sing}=\braces{0}$ et $k_{l}(0)>k$, avec $k_{l}(0)=ord_{0}(\psi_{l})$ ($ord_{0}(\psi_{l}$)
 est l'ordre en 0 de $\psi_{l},\;l\in\braces{2,\ldots,n}$). Nous allons commencer par montrer l'estimation suivante
   $$
     \norm{\varphi(z)}_{2}\pgoe c_{\varphi}\modul{z}^{k},\quad \forall z \in D(0,\rho),
   $$
   où $c_{\varphi}$ est une constante réelle strictement positive ne dépendant que de $\varphi$, la constante
   $\rho$ est fixée dans $\rbrack 0,1 \lbrack$ dépendant de $\varphi$. On peut supposer que $\varphi$ est
   holomorphe dans $\overline{D(0,1)}$, on peut donc développer les fonctions composantes $\varphi_{l}$
    en série entière dans $\overline{D(0,1)}$ en $0$,
   $$
     \psi_{l}(z)=\sum_{j=k_{l}(0)}^{\infty}\frac{1}{j!}\frac{\partial^{j}}{\partial
     z^{j}}\psi_{l}(0)z^{k},
   $$
   on obtient en factorisant par $z^{k_{l}(0)}$ et les inégalités triangulaires:
     $$
       \modul{\psi_{l}(z)}\pgoe\modul{z^{k_{l}(0)}}\parentheses{\frac{1}{k_{l}(0)!}\modul{\frac{\partial^{k_{l}(0)}
       }{\partial
     z^{k_{l}(0)}}\psi_{l}(0)}-\modul{\sum_{j=k_{l}(0)+1}^{\infty}\frac{\partial^{j}
       }{\partial
     z^{j}}\psi_{l}(0)z^{j}}},
     $$
     $\forall z \in D(0,1)$. Avec les inégalités de Cauchy:
     $$
\modul{\sum_{j=k_{l}(0)+1}^{\infty}\frac{\partial^{j}
       }{\partial
     z^{j}}\psi_{l}(0)z^{j}}\ppoe\sum_{j=k_{l}(0)+1}^{\infty}\rho^{j-k_{l}(0)}\norm{\psi_{l}}_{D(0,1)},\;\forall
     z\in D(0,\rho).
     $$
     La constante $\rho$ est dans l'intervalle $\rbrack 0,1\lbrack$. On a donc l'estimation
     ci-dessous
     $$
\modul{\sum_{j=k_{l}(0)+1}^{\infty}\frac{\partial^{j}
       }{\partial
     z^{j}}\psi_{l}(0)z^{j}}\ppoe\norm{\psi_{l}}_{D(0,1)}\frac{\rho}{1-\rho}.
     $$
     Si maintenant, nous faisons tendre $\rho\rightarrow 0$ le membre de gauche
     tend à son tour vers $0$, il existe donc $\rho_{l}\in \rbrack 0,1\lbrack$, tel que
     $$
\modul{\sum_{j=k_{l}(0)+1}^{\infty}\frac{\partial^{j}
       }{\partial
     z^{j}}\psi_{l}(0)z^{j}}\ppoe \frac{1}{2k_{l}(0)!}\modul{\frac{\partial^{k_{l}(0)}
       }{\partial
     z^{k_{l}(0)}}\psi_{l}(0)},\quad\forall z\in D(0,\rho_{l}).
     $$
     Nous obtenons l'inégalité suivante
     $$
       \modul{\psi_{l}(z)}\pgoe\modul{z}^{k_{l}(0)}\frac{1}{2k_{l}(0)!}\modul{\frac{\partial^{k_{l}(0)}
       }{\partial
     z^{k_{l}(0)}}\psi_{l}(0)},\quad\forall l \in\braces{1,\ldots,n},\;\forall z\in D(0,\rho),
     $$
où $\rho=\min_{1\ppoe l \ppoe n}(\rho_{l})$ assez petit. Il s'en suit
     \begin{equation}
       \sqrt{\modul{cz^k}^2+\modul{\psi_{2}(z)}^2+\cdots+\modul{\psi_{n}(z)}^2}
       \pgoe c\modul{z}^{k},\quad\forall z \in D(0,\rho).
     \end{equation}
      Fixons $\hat{z}_{1}$ dans $D(0,1)$.
     Soit maintenant, $\gamma$ un chemin de $CVB^{(0,\varphi(\hat{z}_{1}))}(\lbrack
     0,1\rbrack,\tilde{A})$. Choisissons $\sigma=(t_{0},\ldots,t_{p})$ une subdivision dans
     $\mathfrak{S}_{[0,1]}$. Par construction de la paramétrisation de Puiseux, voire le
     paragraphe 3.4 et en particulier la proposition 2, il existe $z_{0},\ldots,z_{p}$ dans $D(0,1)$,
      tels que
     $\varphi(z_{j})=\gamma(t_{j})$, pour tout $j$ dans $\braces{0,\ldots,p}$ avec
     $\modul{z_{j}}<\modul{z_{j+1}}$ et $z_{p}=\hat{z}_{1}$.
     $$
        V_{\sigma_{}}(\gamma)=\sum_{j=0}^{j=p-1}\norm{\varphi(z_{j+1})-\varphi(z_{j})}_{2}.
     $$
     On peut trouver une subdivision $\sigma_{0}$ de
      $\mathfrak{S}_{[0,1]}$ plus fine que $\sigma$ de
     sorte que $z_{1}$ soit dans $D(0,\rho)$.
     Donc d'après l'estimation (11) on a
     \begin{equation}
       V_{\sigma_{0}}(\gamma)\pgoe
       c\modul{z_{1}}^{k}+\sum_{j=1}^{j=p-1}\norm{\varphi(z_{j+1})-\varphi(z_{j})}_{2}.
     \end{equation}
     Rappelons que $z_{0}=0$. Il nous reste à minorer les termes de la somme. Nous avons
     $$
        \varphi(z_{j+1})-\varphi(z_{j})=\int_{0}^{1}\frac{d}{dt}\varphi(tz_{j+1}+(1-t)z_{j})dt,\quad\forall
        j\in\braces{1,\ldots,p-1},
     $$
     d'où
     $$
        \norm{\varphi(z_{j+1})-\varphi(z_{j})}_{2}=\modul{z_{j+1}-z_{j}}\parentheses{\sum_{l=1}^{l=n}
        \modul{\int_{0}^{1}\prim{\varphi_{l}}(tz_{j+1}+(1-t)z_{j})dt}^{2}}^{\frac{1}{2}},
     $$
     où les $\varphi_l$ sont les fonctions composantes de $\varphi$,
donc
     $$
        \norm{\varphi(z_{j+1})-\varphi(z_{j})}_{2}\pgoe\modul{z_{j+1}-z_{j}}
        \modul{\int_{0}^{1}k(tz_{j+1}+(1-t)z_{j})^{k-1}dt}.
     $$
En conclusion,
     $$
       \norm{\varphi(z_{j+1})-\varphi(z_{j})}_{2}\pgoe\modul{(z_{j+1})^{k}-(z_{j})^{k}}
       \pgoe\modul{z_{j+1}}^{k}-\modul{z_{j}}^{k},\quad\forall j\in\braces{1,\ldots,p-1},
     $$
     car $\modul{z_{j}}<\modul{z_{j+1}},\;\forall j \in \braces{0,\ldots,p-1}$.
En reportant l'inégalité ci-dessus dans l'estimation (12), on obtient
     $$
       V_{\sigma_{0}}(\gamma)\pgoe \min(c,1)\modul{z_{p}}^{k}=\min(c,1)\modul{\hat{z}_{1}}^{k}.
     $$
L'inégalité ci-dessus est encore vraie pour tout
$\sigma\in\mathfrak{S}_{[0,1]},\;\sigma\succcurlyeq\sigma_{0}$; d'où
     $$
         V(\gamma)\pgoe\min(c,1)\modul{\hat{z}_{1}}^{k}.
     $$
Nous pouvons conclure
     $$
       d(0,\varphi(\hat{z}_{1}))\pgoe c_{1}\modul{\hat{z}_{1}}^{k},\quad c_1=\min(c,1).
     $$
     En ce qui concerne l'autre inégalité, il suffit de remarquer que nous avons,
     $$
       d(0,\varphi(\hat{z}_{1}))\ppoe\int_{0}^{1}\norm{\frac{d}{dt}\varphi(t\hat{z}_{1})}_{2}dt,
     $$
     et ceci par la construction même de la métrique $d$.
     Par un calcul similaire aux précédents, nous déduisons que
     $$
       \int_{0}^{1}\norm{\frac{d}{dt}\varphi(t\hat{z}_{1})}_{2}dt\ppoe
       c_{2}\modul{\hat{z}_{1}}^{k},
     $$
où $c_{2}$ est une constante ne dépendant que de $\varphi$. $\blacksquare$

Nous allons démontrer maintenant que la fonction extrémale $V_{A\cap
B(x_{0},c_{1}\eps^{k})}$ vérifie la propriété d'$HCP$ dans $\tilde{A}$ muni de la
métrique des géodésiques définie dans ce paragraphe.
       \begin{pro}
             Soit $A$ une courbe algébrique réelle localement irréductible muni de $d(\cdot,\cdot)$ la métrique des
             géodésiques dans $\tilde{A}$. Alors, la fonction extrémale
              vérifie la propriété de continuité de Hölder locale dans $\tilde{A}$, pour la métrique
             $d(\cdot,\cdot)$. Plus précisément, il existe $\eps_{0},\,\mu_0\in]0,1[$ (dépendant de $x_0$)
             tel que $\forall\eps\in]0,\eps_{0}[,\;\forall\mu\in\obracks{0,\mu_0}$ :

$\bullet Si\;x_{0}\in A_{sing},$

       \begin{equation}
              V_{A\cap B(x_{0},c_{1}\eps^{k})}(\xi)\ppoe C(x_{0})\mu^{\frac{1}{2k}},
             \quad\forall \xi\in B_{\tilde{A}}(x_{0},\eps^{k}\mu)
       \end{equation}
où $k$ est la multiplicité complexe du point singulier $x_{0}$ dans $\tilde{A}$.

$\bullet Si\;x_{0}\in A_{reg}\cap\partial A,$

              \begin{equation}
              V_{A\cap
             B(x_{0},c_{1}\eps)}(\xi)\ppoe C(x_{0})\mu^{\frac{1}{2}},
             \quad \forall \xi\in B_{\tilde{A}}(x_{0},\eps\mu).
       \end{equation}

$\bullet Si\; x_{0}\in A_{reg}\setminus\partial A,$
       \begin{equation}
             V_{A\cap B(x_{0},c_{1}\eps)}(\xi)\ppoe
            C(x_{0})\mu,\quad
            \forall \xi\in B_{\tilde{A}}(x_{0},\eps\mu).
       \end{equation}

Avec $B_{\tilde{A}}(x_{0},r):=\braces{\xi\in \tilde{A}:d(x_{0},\xi)<r}$, $r$ est un réel
strictement positif et $C(x_{0})$ une constante strictement positive localement
supérieurement majorée.
       \end{pro}

\textbf{Preuve:} Commençons par démontrer $(13)$. Comme nous l'avons fait dans les
démonstrations précédentes, nous pouvons supposer, sans perdre de généralité, que
$x_{0}=0$. D'après la proposition 5 du 4.2, modulo un changement de coordonnées unitaire,
il existe des constantes $\rho$, $\eps_{0}$ et $c_{1}$,
 réelles, strictement positives et dépendant seulement de la paramétrisation locale $\varphi$. Cette
 paramétrisation est la paramétrisation de Puiseux construite dans la proposition 2 du
 3.4; rappelons que nous avons $\rho>\eps_{0}>\eps>0$.

 Soit $\xi\in B_{\tilde{A}}(x_{0},\eps^{k}\mu)$,
 d'après le $(10)$ du lemme 5 du 4.3, l'inégalité suivante
 est vérifiée,
 $\modul{\hat{z}}^{k}\ppoe\frac{\eps^{k}\mu}
 {c_{1}}$, où $\varphi(\hat{z})=\xi$.
 Maintenant, en utilisant le (5) de la proposition 5 du 4.2, on peut écrire:$\forall
 \eps\in]0,\eps_{0}]$, $\forall \theta\in\braces{\frac{2k_{j}^{\sigma}\pi}{k}:j\in\braces{1,\ldots,
 r_{A}^{\sigma}}}$,
    $$
      V_{A\cap B(0,\tilde{c}_{1}\eps^{k})}\circ\varphi(z)\ppoe \tilde{c}_{2}V_{\eps\mathcal{R}^{<r_{A}^{\sigma}>}}(z)\ppoe
      \tilde{c}_{2}V_{\eps I}(e^{-i\theta}z),\quad\forall z\in D(0,\rho)\quad (\square_{1}),
    $$
    où $\sigma\in\braces{+,-}$ et $0<\eps_{0}<\rho$. Choisissons $0<\mu_{0}<c_1$ de sorte que $\forall\mu\in]0,\mu_{0}]$
    implique $\frac{\mu}{c_{1}}<\rho^k$, donc $\parentheses{\frac{\modul{\hat{z}}}{\eps}}^k<\frac{\mu}{c_1}<\rho^k$.
     Avec l'inégalité $(\square_{1})$ et du fait que l'on puisse utiliser le lemme 3 du 4.2, car
     on a choisit $\mu_0$ de sorte que $\frac{\mu_0}{c_1}<1$ nous déduisons:
     $$
       V_{A\cap B(0,\widetilde{c}_{1}\eps^{k})}(\xi)\ppoe \widetilde{c}_2V_{I}\parentheses{e^{-i\theta}\frac{\hat{z}}{\eps}}
       \ppoe C\mu^{\frac{1}{2k}}.
     $$
     Pour démontrer (14) la technique de démonstration est identique, en effet soit $\xi\in
     B_{\tilde{A}}(0,\eps\mu)$, il existe $\hat{z}\in D(0,1)$ tel que $\varphi(\hat{z})=\xi$.
     D'après le (9) du lemme 5 on a $d(0,\xi)\pgoe c_{1}\modul{\hat{z}}$, d'où
     $\modul{\hat{z}}\ppoe\frac{\eps\mu}{c_{1}}$. On conclut avec le lemme 4 du 4.2,
          $$
                   V_{A\cap B(0,\widetilde{c}_1\eps)}(\xi)\ppoe \widetilde{c}_{2}V_{\eps I}(\hat{z})\ppoe C\mu^{\frac{1}{2}}.
          $$
          L'Estimation (15) est semblable à la (13).
      $\blacksquare$
\section{Démonstration des théorèmes}
Nous allons démontrer les théorèmes 1 et 2 qui font l'objet de notre papier. Pour ce
faire, nous commencerons par démontrer le lemme suivant.
\begin{lem}
Soit $\varphi : D(0, r) \longrightarrow \com^{n}$  une application holomorphe définie
dans le disque ouvert $D(0,r)$ de $\com$ ($r>0$) telle que $\varphi (0)=0$. Notons
 $\varphi =\uplet{\varphi}{n}$ et $k= \min_{1\leqslant j \leqslant n}(ord_{0}(\varphi_{j}))$.
 Alors il existe un réel strictement positif $r_{0}$ dépend
 uniquement de l'application $\varphi$, tel que :
 $ \forall p \in \polyn{\ree}{x}{n} , \;$ $ \forall r \in \rbrack 0,r_{0} \lbrack
      , \; \forall z \in D(0,\frac{r_{0}}{2})$,
\begin{equation}
\frac{\modul{\derive{(p \circ \varphi)(z)}{z}}}{\sqrt{\underset{1\ppoe j\ppoe
n}{\sum}\modul{\prim{\varphi_j}(z)}^2}} \leqslant \frac{c}{\pi}\int_{0}^{2  \pi} \frac{d
\theta}{\modul{z + \frac{r}{2}e^{i \theta}}^{k-1}} \times \frac{\norm{p \circ
\varphi}_{D(0,r_0)}}{r}.
\end{equation}
\end{lem}

\textbf{Preuve:}
 Les fonctions composantes de $\varphi$ se développent en série entière car $\varphi$ est
holomorphe, donc :  $\forall l \in \braces{1,\ldots,n}$,
     $$
        \varphi_{l}(\xi)= \sum_{k_{l}\leqslant j}
       a_{l,j}\xi^{j}, \quad k_{l}=ord_{0}(\varphi_{l}), \; k_{l}\geqslant 1
      $$
 Comme $a_{l,k_{l}} \not =
0$ lorsque $z \longrightarrow 0$, il existe $r_{l} \in\rbrack 0,1 \lbrack$ tel que :
      $$
       \modul{
       \varphi_{l}(\xi)} \geqslant c_{\varphi ,l} \modul{\xi}^{k_{l}-1},\quad \forall \xi \in D(0,
       r_{l}),
      $$
      d'où
      $$
          \norm{\prim{\varphi}(\xi)}_{2}\geqslant c_{\varphi} \modul{\xi}^{k-1}, \quad \forall \xi
           \in D(0,r_{0}),\quad r_{0}=\min _{1\leqslant l\leqslant n} ({r_{l}}), \quad k=\min_{1 \leqslant  j
          \leqslant n}({k_{l}}).
      $$
      Majorons $ \modul{\frac{1}{z^{k-1}}\derive{(p \circ \ \varphi)}{z}}$. Il
est aisé de voir que la singularité en $0$ est artificielle. Donc
$\frac{1}{z^{k-1}}\derive{(p \circ  \varphi)}{z}$  est holomorphe dans $D(0,r_{0})$. Nous
pouvons appliquer successivement  la formule intégrale de Cauchy à
$\frac{1}{z^{k-1}}\derive{(p \circ  \varphi)}{z}$, afin d'obtenir les inégalités voulues.
     $$
            \frac{1}{z^{k-1}}\derive{(p \circ \ \varphi}{z})(z) =\frac{1}{2  \pi i}
             \int_{C(z,\frac{r}{2})}\frac{1}{\zeta^{k-1}}\frac{1}{2  \pi i} \int_{C(\zeta,\frac{r}{2})}
             \frac{p \circ \varphi (\xi)}{(\zeta - \xi)^{2}} d \xi \frac{d \zeta}{\zeta - z}
     $$
     $$
          =\frac{1}{2
          \pi i} \int_{0}^{2  \pi} \frac{1}{(z + \frac{r}{2}e^{i \theta_{2} })^{k-1}} \frac{1}{2  \pi i}
           \int_{0}^{2  \pi} \frac{p \circ \varphi(z + \frac{r}{2}e^{i \theta_{2} }+\frac{r}{2}e^{i
           \theta_{1} })}{\frac{r}{2}e^{i \theta_{1} }} d  \theta_{1} d  \theta_{2},
     $$
     $$
          \forall r \in
          \rbrack 0,r_{0}\lbrack , \quad \forall z \in D(0, \frac{r}{2}).
     $$
     D'où
     $$
           \modul{\frac{1}{z^{k-1}}\derive{p \circ \varphi}{z} (z)} \leqslant \frac{c}{\pi} \int_{0}^{2
           \pi} \frac{1}{\modul{z + \frac{r}{2}e^{i\theta}}^{k-1}} \frac{\norm{p \circ
           \varphi}_{D(z,r_{0})}}{r}d  \theta,
     $$
     $$
            \forall r \in \rbrack 0,r_{0} \lbrack ,\quad \forall z
            \in D(0, \frac{r}{2}).
     $$
     D'où l'estimation (16). $\blacksquare$

     Commençons par démontrer $1.$ du théorème 1.

     \textbf{Preuve:}
     Sans perdre de généralité, nous pouvons supposer que $x_{0}=0$ et que $0\in A_{sing}$.

     Soit $p$ un polynôme dans $\polyn{\ree}{x}{n}$ fixé.
     Si $v\in C(A,0)$ est un vecteur tangent unitaire, alors la branche de $A$ tangente à $v$
     est paramétriseé par $\restric{\varphi}{\mathcal{R}^{<r_A^+>}}$ ou
     $\restric{\varphi}{\mathcal{R}^{<r_A^->}}$ ($\varphi$ paramétrisation de Puiseux),
      ceci d'après la proposition 2 du 3.4. On peut donc choisir $l$ un entier dans $\braces{0,\ldots,k-1}$, tel que le segment
     $\lbrack 0,e^{i\frac{2l\pi}{k}}\rbrack\subset\mathcal{R}^{<r_A^+>}\cup\mathcal{R}^{<r_A^->}$
     et $\restric{\varphi}{\lbrack 0,e^{i\frac{2l\pi}{k}}\rbrack}$ paramétrise la branche de $A$
      à laquelle le vecteur $v$ est tangent. Pour des commodités d'écriture et de calculs, nous ne changerons rien au résultat voulu, si nous
     transformons le segment $\lbrack 0,e^{i\frac{2l\pi}{k}}\rbrack$ en $[0,1]$, par une rotation dans
     $\com$ d'angle $-\frac{2l\pi}{k}$ (bien évidemment $\com$ est orienté dans le sens direct).
     L'ensemble $\mathcal{R}^{<r_A^+>}\cup\mathcal{R}^{<r_A^->}$ n'est en rien modifié par cette rotation.

     D'après l'inégalité (16) lemme 6:
     \begin{equation}
                      \modul{D_{v}p(0)}\ppoe
                      \frac{c}{r^{k}}\norm{p\circ\varphi}_{D(0,r)},\quad\forall
                      r\in]0,r_{0}],
     \end{equation}
la constante $r_{0}$ est réelle strictement positive ne dépendant que de $\varphi$.
D'après le (13) la proposition 6, il existe $\mu_0,\eps_0>0$ tel que: $\forall
\xi\in\tilde{A},\;d(\xi,0)\ppoe\eps^k\mu$,
$$
           V_{A\cap B(x_0,c_1\eps^k)}(\xi)\ppoe C(x_0)\mu^{\frac{1}{2k}}\quad (\square_1)
$$
 Quitte à diminuer la valeur de $\mu_{0}$, on peut donc écrire:
    $$
      \modul{D_{v}p(0)}\ppoe
      \frac{c}{(\eps\mu)^{k}}\norm{p\circ\varphi}_{D(0,\eps\mu)},\quad\forall \mu\in
      ]0,\mu_{0}[,\;\forall \eps\in ]0,\eps_{0}[.
    $$
     Estimons $\norm{p\circ\varphi}_{D(0,\eps\mu)}$; c'est du (10) du lemme 5 du 4.3
    dont nous obtenons l'estimation ci-dessous,
    \begin{equation}
       \norm{p\circ\varphi}_{D(0,\eps\mu)}\ppoe \norm{p}_{B_{\tilde{A}}(0,c_{2}(\eps\mu)^{k})},\quad\forall \mu\in
      ]0,\mu_{0}[,\;\forall \eps\in ]0,\eps_{0}[,
    \end{equation}
    (car $\varphi(D(0,\eps\mu))\subset B_{\tilde{A}}(0,c_{2}(\eps\mu)^{k})$).
    $$
     \modul{D_{v}p(0)}\ppoe
      \frac{c}{(\eps\mu)^{k}}\norm{p}_{B_{\tilde{A}}(0,c_{2}(\eps\mu)^{k})},\quad\forall \mu\in
      ]0,\mu_{0}[,\;\forall \eps\in ]0,\eps_{0}[,
    $$
avec inégalité de Bernstein-Walsh:
    $$
   \modul{D_{v}p(0)}\ppoe \frac{c}{(\eps\mu)^k}\norm{p}_{A\cap
    B(0,c_{1}\eps^{k})}\exp\parentheses{\deg(p)\sup_{B_{\tilde{A}}(0,c_{2}(\eps\mu)^{k})}V_{A\cap
    B(0,c_{1}\eps^{k})}}.
    $$
    En utilisant la propriété d'$HCP$ de la fonction de Green ($\square_1$), on obtient:
    $\forall \mu\in ]0,\mu_{0}[,\;\forall \eps\in ]0,\eps_{0}[,$
    $$
        \modul{D_{v}p(0)}\ppoe \frac{c}{(\eps\mu)^{k}}\norm{p}_{A\cap
        B(0,c_{1}\eps^{k})}\exp\parentheses{C_3\deg(p){(\mu^{k})^{\frac{1}{2k}}}},\;
    $$
    d'où
    $$
        \modul{D_{v}p(0)}\ppoe \frac{c}{(\eps\mu)^{k}}\norm{p}_{A\cap
        B(0,c_{1}\eps^{k})}\exp\parentheses{C_3\deg(p){\mu^{\frac{1}{2}}}}.
    $$
En posant $\mu=\frac{\tilde{C}}{(\deg p)^2}$, avec $0<\tilde{C}<\mu_0$
    $$
         \modul{D_{v}p(0)}\ppoe c\,e^{C_3\cdot\sqrt{\tilde{C}}}\parentheses{\frac{\tilde{C}(\deg(p))^{2}}{\eps}}^{k}\norm{p}_{A\cap
        B(0,c_{1}\eps^{k})}.
    $$
Le (1) du théorème 1 est donc montré.

    Pour démontrer $2.$ et $3.$ du théorème 1, la  démonstration techniquement identique, seulement il faut utiliser le $(14)$ et $(15)$
    de la proposition 6, le (9) et (20) du lemme 5 et le lemme 4. $\blacksquare$

    Un exemple intéressant a été démontré par L.Bos \cite{bos3} d'inégalités de
Markov tangentielles sur certaines courbes algébriques de $\ree^{2}$. Il montre que
l'exposant $k$ est optimum sur les courbes de type:
   $$
      \Gamma=\braces{(t^{p},t^{q}):t\in [0,1]},
   $$
où p et q sont deux entiers naturels premiers entre eux, tels que $p<q$. Pour des
inégalités globales, il montre que $k$ ne peut être plus petit que $p$. Or cet entier $k$
est aussi la multiplicité complexe du point $(0,0)$ de la courbe algébrique complexifiée
$\tilde{\Gamma}$ de $\Gamma$, sans oublier que tous les autres points de $\Gamma$ sont
des points réguliers, donc localement ont un exposant de Markov au moins égal à 1 ou $p$.
Le théorème 1 conforte l'idée que la multiplicité des points d'une courbe algébrique joue
un rôle sur l'exposant de Markov, on peut même penser que les points singuliers d'une
courbe algébrique réelle se comportent comme un bord, ce qui expliquerait l'apparition du
$2k$ à l'exposant.
\\

L'auteur remercie chaleureusement M. Baran, J. Duval, L. Bos, N. Levenberg et A. Zeriahi
pour les entretiens fructueux et intéressants qui ont permis d'achever ce travail.

\nocite{*}
\bibliographystyle{plain}
\bibliography{Markov.bib}
Adresse:\\
Laboratoire E.Picard,\\ U.M.R. C.N.R.S. 5580,\\
Département de Mathématiques,\\
Université Paul Sabatier\\
118, route de Narbonne 31062\\
TOULOUSE CEDEX 04\\
FRANCE.\\
 Mail:\url{gendre@picard.ups-tlse.fr}

\end{document}